# CONVERGENCE OF ALGORITHMS FOR RECONSTRUCTING CONVEX BODIES AND DIRECTIONAL MEASURES[1]

BY RICHARD J. GARDNER, MARKUS KIDERLEN
AND PEYMAN MILANFAR

*Western Washington University, University of Aarhus
and University of California, Santa Cruz*

We investigate algorithms for reconstructing a convex body $K$ in $\mathbb{R}^n$ from noisy measurements of its support function or its brightness function in $k$ directions $u_1, \ldots, u_k$. The key idea of these algorithms is to construct a convex polytope $P_k$ whose support function (or brightness function) best approximates the given measurements in the directions $u_1, \ldots, u_k$ (in the least squares sense). The measurement errors are assumed to be stochastically independent and Gaussian.

It is shown that this procedure is (strongly) consistent, meaning that, almost surely, $P_k$ tends to $K$ in the Hausdorff metric as $k \to \infty$. Here some mild assumptions on the sequence $(u_i)$ of directions are needed. Using results from the theory of empirical processes, estimates of rates of convergence are derived, which are first obtained in the $L_2$ metric and then transferred to the Hausdorff metric. Along the way, a new estimate is obtained for the metric entropy of the class of origin-symmetric zonoids contained in the unit ball.

Similar results are obtained for the convergence of an algorithm that reconstructs an approximating measure to the directional measure of a stationary fiber process from noisy measurements of its rose of intersections in $k$ directions $u_1, \ldots, u_k$. Here the Dudley and Prohorov metrics are used. The methods are linked to those employed for the support and brightness function algorithms via the fact that the rose of intersections is the support function of a projection body.

Received October 2004; revised July 2005.
[1]Supported in part by NSF Grants DMS-02-03527 and CCR-99-84246, and by the Carlsberg Foundation.

**1. Introduction.** The problem of reconstructing an unknown shape from a finite number of noisy measurements of its support function [giving the (signed) distances from some fixed reference point, usually taken to be the origin, to the support hyperplanes of the shape] has attracted much attention. The nature of the measurements makes it natural to restrict attention to convex bodies. Prince and Willsky [27] used such data in computerized tomography as a prior to improve performance, particularly when only limited data is available. Lele, Kulkarni and Willsky [21] applied support function measurements to target reconstruction from range-resolved and Doppler-resolved laser-radar data. The general approach in these papers is to fit a polygon or polyhedron to the data by a least squares procedure. In contrast, Fisher, Hall, Turlach and Watson [8] use spline interpolation and the so-called von Mises kernel to fit a smooth curve to the data in the two-dimensional case. This method was taken up by Hall and Turlach [16] and Mammen, Marron, Turlach and Wand [22], the former dealing with convex bodies with corners and the latter giving an example to show that the algorithm of Fisher, Hall, Turlach and Watson [8] may fail for a given data set. Further applications and the three-dimensional case can be found in papers by Gregor and Rannou [14], Ikehata and Ohe [18] and Karl, Kulkarni, Verghese and Willsky [19].

Despite all this work, the convergence of even the most straightforward of the reconstruction algorithms has apparently never been proved. In Theorem 6.1 below, we provide such a proof for an algorithm we call Algorithm NoisySupportLSQ, due to Prince and Willsky [27]. By convergence, we mean that, given a suitable sequence of directions, the estimators, convex polytopes, obtained by running the algorithm with noisy measurements taken in the first $k$ directions in the sequence, converge in suitable metrics (the $L_2$ and Hausdorff metrics) to the unknown convex body as $k$ tends to infinity. Suitable sequences of directions are those that are "evenly spread," only slightly more restrictive than the obviously necessary condition that the sequence be dense in the unit sphere.

Moreover, by applying some beautiful and deep results from the theory of empirical processes, we are able to provide in Theorem 6.2 estimates of rates of convergence of the estimators to the unknown convex body. Some considerable technicalities are involved, and some extra conditions are required, of which, however, only a rather stronger condition on the sequence of directions should be regarded as really essential. Convergence rates depend on the dimension of the unknown convex body; for example, for the $L_2$ metric, the rate is of order $k^{-2/5}$ in the two-dimensional case, and $k^{-1/3}$ in the three-dimensional case.

Analogous results are obtained for an algorithm, Algorithm NoisyBrightLSQ, essentially that proposed recently by Gardner and Milanfar [13], that constructs an approximating convex polytope to an unknown origin-symmetric convex body from a finite number of noisy measurements of its brightness function (giving the areas of the shadows of the body on hyperplanes). The very existence of such an algorithm is highly nontrivial, due to the extremely weak data; each measurement is a single scalar that provides no information at all about the shape of the shadow! Nevertheless, the algorithm has been successfully implemented, even in three dimensions. Here we are able to prove, for the first time, convergence (Theorem 7.2), with estimates of rates of convergence (Theorem 7.6) also for this algorithm. One technical device that aids in this endeavor is the so-called projection body, whose support function equals the brightness function of a given convex body. This allows some of our results on reconstruction from support functions to be transferred to the new reconstruction problem. However, we require additional deep results on projection bodies (a subclass of the class of zonoids) from the theory of convex geometry due to Bourgain and Lindenstrauss [1] and Campi [4]. Examples of rates of convergence we obtain, for the Hausdorff metric, are of order $k^{-4/15}$ in the two-dimensional case and $k^{-1/30}$ in the three-dimensional case.

Most of our results are actually much more informative in that they indicate also how the convergence depends on the noise level and the scaling of the input body. A discussion and the results of some Monte Carlo simulations can be found in Section 8.

Many auxiliary results are obtained in the course of proving the convergence of these algorithms, but one is perhaps worth special mention. Roughly speaking, the results we employ from the theory of empirical processes give rates of convergence of least squares estimators to an unknown function in terms of the metric entropy of the class of functions involved. In obtaining our results on reconstruction from support functions, it turns out that we therefore need an estimate of the metric entropy of the class of compact convex subsets of the unit ball $B$ in $n$-dimensional space, with the Hausdorff metric. Luckily, the precise order of this, $t^{-(n-1)/2}$ for sufficiently small $t > 0$, was previously established by Bronshtein [3] (see Proposition 5.4; it is traditional to talk of $\varepsilon$-entropy rather than $t$-entropy, but we require $\varepsilon$ for a different purpose in this paper). In the problem of reconstruction from brightness functions, however, we need to know the metric entropy of the class of origin-symmetric zonoids contained in $B$. As far as we know, this natural problem has not been addressed before. For $n = 2$, it is easy to see that the answer, $t^{-1/2}$, is unchanged, but in Theorem 7.3 we show that, for fixed $n \geq 3$ and any $\eta > 0$, the $t$-entropy of this class is $O(t^{-2(n-1)/(n+2)-\eta})$ for sufficiently small $t > 0$. This is somewhat remarkable, since the $t$-entropy becomes $O(t^{-2})$ as $n$ tends to infinity, in complete contrast to the case of general compact convex sets. The hard work behind Theorem 7.3 is done in the highly technical papers of Bourgain and Lindenstrauss [2] and Matoušek [24] on the approximation of zonoids by zonotopes.

While most of the paper is devoted to reconstruction of convex bodies, Section 9 focuses on a problem from stereology, that of reconstructing an unknown directional measure of a stationary fiber process from a finite number of noisy measurements of its rose of intersections. It turns out that the corresponding algorithm, Algorithm NoisyRoseLSQ, is very closely related to Algorithm NoisyBrightLSQ, due to the fact that the rose of intersections is the support function of a projection body. This fact was also used by Kiderlen [20], where an estimation method similar to Algorithm NoisyRoseLSQ was suggested and analyzed. Convergence of Algorithm NoisyRoseLSQ was proved by Männle [23], but also follows easily from our earlier results (see Proposition 9.1). With suitable extra assumptions, we can once again obtain estimates of rates of convergence of the approximating measures to the unknown directional measure. These are first given for the Dudley metric in Theorem 9.4, but can easily be converted to rates for the Prohorov metric. For example, for the Prohorov metric, the rate is of order $k^{-1/20}$ in the three-dimensional case.

**2. Definitions, notation and preliminaries.** As usual, $S^{n-1}$ denotes the unit sphere, $B$ the unit ball, $o$ the origin and $\|\cdot\|$ the norm in Euclidean $n$-space $\mathbb{R}^n$. It is assumed throughout that $n \geq 2$. A *direction* is a unit vector, that is, an element of $S^{n-1}$. If $u$ is a direction, then $u^\perp$ is the $(n-1)$-dimensional subspace orthogonal to $u$. If $x, y \in \mathbb{R}^n$, then $x \cdot y$ is the inner product of $x$ and $y$ and $[x, y]$ denotes the line segment with endpoints $x$ and $y$.

If $A$ is a set, $\dim A$ is its *dimension*, that is, the dimension of its affine hull, and $\partial A$ is its *boundary*. The notation for the usual (orthogonal) *projection* of $A$ on a subspace $S$ is $A|S$. A set is *origin symmetric* if it is centrally symmetric, with center at the origin.

We write $V_k$ for $k$-dimensional Lebesgue measure in $\mathbb{R}^n$, where $k = 1, \ldots, n$, and where we identify $V_k$ with $k$-dimensional Hausdorff measure. If $K$ is a $k$-dimensional convex subset of $\mathbb{R}^n$, then $V(K)$ is its *volume* $V_k(K)$. Define $\kappa_n = V(B)$. The notation $dz$ will always mean $dV_k(z)$ for the appropriate $k = 1, \ldots, n$.

Let $\mathcal{K}^n$ be the family of nonempty compact convex subsets of $\mathbb{R}^n$. A set $K \in \mathcal{K}^n$ is called a *convex body* if its interior is nonempty. If $K \in \mathcal{K}^n$, then

$$h_K(x) = \max\{x \cdot y : y \in K\},$$

for $x \in \mathbb{R}^n$, is its *support function* and

$$b_K(u) = V(K|u^\perp),$$

for $u \in S^{n-1}$, its *brightness function*. Any $K \in \mathcal{K}^n$ is uniquely determined by its support function. If $K$ is an origin-symmetric convex body, it is also uniquely determined by its brightness function. The *Hausdorff distance* $\delta(K,L)$ between two sets $K, L \in \mathcal{K}^n$ can be conveniently defined by

$$\delta(K,L) = \|h_K - h_L\|_\infty.$$

We shall also employ the $L_2$ distance $\delta_2(K,L)$ defined by

$$\delta_2(K,L) = \|h_K - h_L\|_2.$$

By Proposition 2.3.1 of [15], there is a constant $c = c(n)$ such that if $K$ and $L$ are contained in $RB$ for some $R > 0$, then

(1) $$\delta(K,L) \leq cR^{(n-1)/(n+1)} \delta_2(K,L)^{2/(n+1)},$$

which shows (together with a trivial inequality in the converse direction) that both metrics induce the same topology on $\mathcal{K}^n$.

A *zonotope* is a vector sum of finitely many line segments. A *zonoid* is the limit in the Hausdorff metric of zonotopes. The *projection body* of a convex body $K$ in $\mathbb{R}^n$ is the origin-symmetric convex body $\Pi K$ defined by

$$h_{\Pi K} = b_K.$$

An introduction to the theory of projection bodies is provided by Gardner [10], Chapter 4. It turns out that projection bodies are precisely the $n$-dimensional origin-symmetric zonoids. For this reason, we shall denote the set of projection bodies in $\mathbb{R}^n$ by $\mathcal{Z}^n$.

The *surface area measure* $S(K,\cdot)$ of a convex body $K$ is defined for Borel subsets $E$ of $S^{n-1}$ by

(2) $$S(K,E) = V_{n-1}(g^{-1}(K,E)),$$

where $g^{-1}(K,E)$ is the set of points in $\partial K$ at which there is an outer unit normal vector in $E$. The convex body $P$ is a zonotope if and only if $P = \Pi K$ for some origin-symmetric convex polytope $K$. In this case, $S(K,\cdot)$ is a sum of point masses, each located at one of the directions of the line segments whose sum is $P$ and with weight equal to half the length of this line segment. This fact will be used in a reconstruction algorithm in Section 7.

A fundamental result is *Minkowski's existence theorem* (see, e.g., [10], Theorem A.3.2), which says that a finite Borel measure $\mu$ in $S^{n-1}$ is the surface area measure of some convex body $K$ in $\mathbb{R}^n$, unique up to translation, if and only if $\mu$ is not concentrated on any great sphere and

$$\int_{S^{n-1}} u \, d\mu(u) = o.$$

The treatise of Schneider [28] is an excellent general reference for all of these topics.

Let $U = \{u_1, \ldots, u_k\} \subset S^{n-1}$. The *nodes* corresponding to $U$ are defined as follows. The hyperplanes $u_i^\perp$, $i = 1, \ldots, k$, partition $\mathbb{R}^n$ into a finite set of polyhedral cones, which intersect $S^{n-1}$ in a finite set of spherically convex regions. The nodes $\pm v_j \in S^{n-1}$, $j = 1, \ldots, l$, are the vertices of these regions. Thus, when $n = 2$, the nodes are simply the $2k$ unit vectors each of which is orthogonal to some $u_i$, $i = 1, \ldots, k$. When $n = 3$, each $v_j$ is of the form $(u_i \times u_{i'})/\|u_i \times u_{i'}\|$, where $1 \leq i < i' \leq k$ and $u_i \neq \pm u_{i'}$. Thus, for $n = 3$, $l \leq k(k-1)/2$ and in general, $l = O(k^{n-1})$. Campi, Colesanti and Gronchi [5] proved the following result.

PROPOSITION 2.1. *Let $K$ be a convex body in $\mathbb{R}^n$ and let $U = \{u_1, \ldots, u_k\} \subset S^{n-1}$ span $\mathbb{R}^n$. Among all convex bodies with the same brightness function values as $K$ in the directions in $U$, there is a unique origin-symmetric convex polytope $P$, of maximal volume and with each of its facets orthogonal to one of the nodes corresponding to $U$.*

This implies that, for any projection body $\Pi K$ and any finite set of directions $U \subset S^{n-1}$, there is a zonotope $Z$ with $h_Z(u) = h_{\Pi K}(u)$, for all $u \in U$. Moreover, $Z$ can be written as a sum of line segments, each parallel to some node corresponding to $U$.

The following deep result was proved independently by Campi [4] (for $n = 3$) and Bourgain and Lindenstrauss [1]. The latter authors state their theorem in terms of a metric other than the Hausdorff metric, and make an additional assumption on the distance between the projection bodies. Groemer ([15], Theorem 5.5.7) presents the version below, and his proof yields the estimate of the constant in (4). This involves some tedious calculations (see www.ac.wwu.edu/~gardner; no attempt was made to obtain the optimal expression). In (4) and throughout the paper, the "big $O$" notation is used in the sense of "less than a constant multiple depending only on $n$."

PROPOSITION 2.2. *Let $K$ and $L$ be origin-symmetric convex bodies in $\mathbb{R}^n$, $n \geq 3$, such that*

$$r_0 B \subset K, L \subset R_0 B,$$

*for some $0 < r_0 < R_0$. If $0 < a < 2/(n(n+4))$, there is a constant $c = c(a, n, r_0, R_0)$ such that*

(3) $$\delta(K, L) \leq c \delta_2(\Pi K, \Pi L)^a.$$

*Moreover, if $0 < a < 2/(n(n+4))$ is fixed, $r_0 < 1$ and $R_0 > 1$, then*

(4) $$c = O(r_0^{-2n-1} R_0^5).$$

**3. Some properties of sets and sequences of unit vectors.** In this section we gather together some basic results on sets and sequences of unit vectors that will be useful in Sections 5 and 7.

If $\{u_1, \ldots, u_k\}$ is a finite subset of $S^{n-1}$, its *spread* $\Delta_k$ is defined by

$$\Delta_k = \max_{u \in S^{n-1}} \min_{1 \leq i \leq k} \|u - u_i\|. \tag{5}$$

For $i = 1, \ldots, k$, let $\Omega_i$ be the spherical Voronoi cell

$$\Omega_i = \{u \in S^{n-1} : \|u - u_i\| \leq \|u - u_j\| \text{ for all } 1 \leq i, j \leq k\} \tag{6}$$

containing $u_i$. Then $\bigcup_{i=1}^k \Omega_i = S^{n-1}$, and we define

$$\omega_k = \max_{1 \leq i \leq k} V_{n-1}(\Omega_i). \tag{7}$$

By the definition of spread, $\{u_1, \ldots, u_k\}$ is a $\Delta_k$-net in $S^{n-1}$, meaning that, for every vector $u$ in $S^{n-1}$, there is an $i \in \{1, \ldots, k\}$ such that $u$ is within a distance $\Delta_k$ of $u_i$. The existence of $\varepsilon$-nets in $S^{n-1}$ with relatively few points is provided by the following well-known result. It can be proved by induction on $n$ in a constructive way; see, for example, [13], Lemma 7.1.

PROPOSITION 3.1. *For each $\varepsilon > 0$ and $n \geq 2$, there is an $\varepsilon$-net in $S^{n-1}$ containing $O(\varepsilon^{1-n})$ points.*

Now let $(u_i)$ be an infinite sequence in $S^{n-1}$. We retain the notation $\Delta_k$ for the spread of the first $k$ points in the sequence, and similarly for $\omega_k$. We need to consider some conditions on $(u_i)$ that are stronger than denseness in $S^{n-1}$. To this end, for $u \in S^{n-1}$ and $0 < t \leq 2$, let

$$C_t(u) = \{v \in S^{n-1} : \|u - v\| < t\}$$

be the open spherical cap with center $u$ and radius $t$. We call $(u_i)$ *evenly spread* if for all $0 < t < 2$, there is a constant $c = c(t) > 0$ and an $N = N(t)$ such that

$$|\{u_1, \ldots, u_k\} \cap C_t(u)| \geq ck, \tag{8}$$

for all $u \in S^{n-1}$ and $k \geq N$.

The following lemma provides relations between various properties of sequences we need later. A discussion of how these properties relate to the well-known concept of a uniformly distributed sequence can be found in the Appendix of [11].

LEMMA 3.2. *Consider the following properties of a sequence $(u_i)$ in $S^{n-1}$:*

(i) $\Delta_k = O(k^{-1/(n-1)})$.

(ii) $\omega_k = O(k^{-1})$ and $(u_i)$ is dense in $S^{n-1}$.
(iii) $(u_i)$ is evenly spread.
(iv) $(u_i)$ is dense in $S^{n-1}$.

Then (i) $\Rightarrow$ (ii) $\Rightarrow$ (iii) $\Rightarrow$ (iv), and there are sequences with property (i).

PROOF. Assume (i), and let $k \in \mathbb{N}$ and $i \in \{1, \ldots, k\}$. Let $\Omega_i$, $1 \leq i \leq k$, be the Voronoi cells corresponding to the set $\{u_1, \ldots, u_k\}$ defined by (6). Note that $\Omega_i \subset C_{\Delta_k}(u_i)$ and hence,

$$V_{n-1}(\Omega_i) \leq V_{n-1}(C_{\Delta_k}(u_i)) \leq V_{n-1}(D_k(u_i)),$$

where $D_k(u_i)$ is the $(n-1)$-dimensional ball in the tangent hyperplane to $S^{n-1}$ at $u_i$, obtained by the inverse spherical projection (with center $o$) of $C_{\Delta_k}(u_i)$. If $\Delta_k < \sqrt{2}$, then $D_k(u_i)$ has center $u_i$ and radius $r_k = \tan(2\arcsin(\Delta_k/2))$. Therefore,

$$\omega_k = \max_{1 \leq i \leq k} V_{n-1}(\Omega_i) \leq r_k^{n-1} \kappa_{n-1} = O(\Delta_k^{n-1}) = O(k^{-1}).$$

Since it is clear that (i) also implies that $(u_i)$ is dense in $S^{n-1}$, (ii) holds.

Suppose that (ii) holds. Fix $0 < t < 2$ and $u \in S^{n-1}$. Cover $S^{n-1}$ with finitely many open caps $C_j = C_{t/6}(v_j)$, $1 \leq j \leq m$. Since $(u_i)$ is dense in $S^{n-1}$, there is an $N = N(t) \in \mathbb{N}$ such that, for $k \geq N$, any of these caps contains at least one point of $\{u_1, \ldots, u_k\}$. The cap $C_{t/3}(u)$ contains at least one $C_j$, and hence a point $u_{i_0}$ with $1 \leq i_0 \leq N$. Note that $N$ does not depend on $u$.

Fix $k \geq N$ and let $\Omega_i$, $1 \leq i \leq k$, be the Voronoi cells corresponding to the set $\{u_1, \ldots, u_k\}$. If $\Omega_i \cap \text{int}\, C_{t/3}(u) \neq \varnothing$, $i \neq i_0$, there must be a point in $C_{t/3}(u)$ closer to $u_i$ than to $u_{i_0}$. This implies $u_i \in C_t(u)$. Consequently,

$$\text{int}\, C_{t/3}(u) \subset \bigcup \{\Omega_i : \Omega_i \cap \text{int}\, C_{t/3}(u) \neq \varnothing\} \subset \bigcup \{\Omega_i : u_i \in C_t(u)\}.$$

Now (ii) implies that there is a $c' = c'(t)$ such that

$$V_{n-1}(C_{t/3}(u)) \leq \sum_{u_i \in C_t(u)} V_{n-1}(\Omega_i)$$
$$\leq \omega_k |\{i : u_i \in C_t(u)\}|$$
$$\leq \frac{c'}{k} |\{u_1, \ldots, u_k\} \cap C_t(u)|.$$

Since the left-hand side of the previous chain of inequalities does not depend on $u$, this yields (iii). That (iii) implies (iv) is clear.

To obtain a sequence with property (i), observe that, by Proposition 3.1, there is a constant $C$ such that, for each $m \in \mathbb{N}$, there is a set $W_m$ of at most $C2^{m(n-1)}$ unit vectors forming a $2^{-m}$-net. Order the elements of each $W_m$ in an arbitrary fashion, and let $(u_i)$ be the sequence obtained by forming one sequence from these finite sequences $W_1$, $W_2$ and so on in that order. Let

$$N_m = C(2^{n-1} + 2^{2(n-1)} + \cdots + 2^{m(n-1)}) = C2^{n-1}\left(\frac{2^{m(n-1)} - 1}{2^{n-1} - 1}\right).$$

Then for all $k \geq N_m$, the points $u_1, \ldots, u_k$ form a $2^{-m}$-net.

Now suppose that $k$ is the least integer such that the points $u_1, \ldots, u_k$ have spread $\Delta_k$, where

$$2^{-m} \leq \Delta_k < 2^{1-m}.$$

Then

$$k \leq N_m = C2^{n-1}\left(\frac{2^{m(n-1)} - 1}{2^{n-1} - 1}\right) < C2^{n-1}\left(\frac{2^{n-1}\Delta_k^{1-n} - 1}{2^{n-1} - 1}\right),$$

or

$$\Delta_k \leq \left(\frac{k(2^{n-1} - 1)}{C2^{2(n-1)}} + \frac{1}{2^{n-1}}\right)^{-1/(n-1)} = O(k^{-1/(n-1)}). \qquad \square$$

Let $(u_i)$ be a sequence of vectors in $S^{n-1}$. For application in Section 7, we need to consider properties of the "symmetrized" sequence

(9) $$(u_i^*) = (u_1, -u_1, u_2, -u_2, \ldots).$$

Let

(10) $$\Delta_k^* = \max_{u \in S^{n-1}} \min_{1 \leq i \leq k} \min\{\|u - u_i\|, \|u - (-u_i)\|\}$$

be the *symmetrized spread* of $u_1, \ldots, u_k$. Also, let $\omega_k^*$ be the maximum $V_{n-1}$-measure of the $2k$ spherical Voronoi cells corresponding to the set $\{\pm u_1, \pm u_2, \ldots, \pm u_k\}$.

Following [20], page 14, we call $(u_i)$ *asymptotically dense* if

$$\liminf_{k \to \infty} \frac{1}{k} |\{u_1, \ldots, u_k\} \cap G| > 0,$$

for all origin-symmetric open sets $G \neq \varnothing$ in $S^{n-1}$.

LEMMA 3.3. *Consider the following properties of a sequence $(u_i)$ in $S^{n-1}$:*

(i) $\Delta_k^* = O(k^{-1/(n-1)})$.
(ii) $\omega_k^* = O(k^{-1})$ and $(u_i^*)$ is dense in $S^{n-1}$.
(iiia) $(u_i^*)$ is evenly spread.
(iiib) $(u_i)$ is asymptotically dense.
(iv) $(u_i^*)$ is dense in $S^{n-1}$.

*Then* (i) $\Rightarrow$ (ii) $\Rightarrow$ (iiia) $\Leftrightarrow$ (iiib) $\Rightarrow$ (iv), *and there are sequences with property* (i).

PROOF. The implications (i) $\Rightarrow$ (ii) $\Rightarrow$ (iiia) $\Rightarrow$ (iv) are direct consequences of Lemma 3.2 and the definition of $(u_i^*)$. The existence statement also follows from this lemma, as any sequence with $\Delta_k = O(k^{-1/(n-1)})$ satisfies $\Delta_k^* = O(k^{-1/(n-1)})$.

That (iiia) implies (iiib) is trivial. To prove the converse, let $C_t(u)$ be an open cap in $S^{n-1}$ of radius $t$, and cover the compact set $S^{n-1}$ with open caps $C_1, \ldots, C_m$ of radius $t/2$. Then $C_j \subset C_t(u)$ for some $j$. If $(u_i)$ is asymptotically dense, we can apply the definition of this property with $G = C_j \cup (-C_j)$ to conclude that there are a constant $c' > 0$ and an $N'$ such that

$$|\{u_1, \ldots, u_k\} \cap (C_j \cup (-C_j))| \geq c'k$$

for all $k \geq N'$ and, hence, that

$$|\{\pm u_1, \ldots, \pm u_k\} \cap C_t(u)| \geq c'k$$

for all $k \geq N'$. From this, it follows easily that $(u_i^*)$ is evenly spread. □

**4. Metric entropy and convergence rates for least squares estimators.** Let $\mathcal{G} \neq \varnothing$ be a class of measurable real-valued functions defined on a subset $E$ of $\mathbb{R}^n$. Suppose that $x_i \in E$, $i = 1, 2, \ldots$, are fixed, and let $X_i$, $i = 1, 2, \ldots$, be independent random variables with mean zero and finite variance. If $g_0 \in \mathcal{G}$, we regard the quantities

$$y_i = g_0(x_i) + X_i,$$

$i = 1, 2 \ldots$, as measurements of the unknown function $g_0$. For $k \in \mathbb{N}$, any function $\hat{g}_k \in \mathcal{G}$ satisfying

$$(11) \qquad \hat{g}_k = \arg\min_{g \in \mathcal{G}} \sum_{i=1}^{k} (y_i - g(x_i))^2$$

is called a *least squares estimator for $g_0$ with respect to $\mathcal{G}$, based on measurements at $x_1, \ldots, x_k$*. (Since $\hat{g}_k$ depends on $y_1, \ldots, y_k$, it also depends on the random variables $X_1, \ldots, X_k$, but this is not made explicit.) If $k$, $\mathcal{G}$ and $x_1, \ldots, x_k$ are clear from the context, we shall simply refer to $\hat{g}_k$ as a least squares estimator for $g_0$. In the definition of $\hat{g}_k$, $x_i$ and $y_i$ are not needed for $i > k$, but later we shall take additional measurements into account in order to examine the asymptotic behavior of $\hat{g}_k$ as $k$ increases. In general, $\hat{g}_k$ need not be unique and the existence of a least squares estimator has to be assumed. In the applications we have in mind, a least squares estimator always exists. To provide the necessary measurability for the background theory to work, a suitable condition can be imposed on the class $\mathcal{G}$. Following [25], page 196, we call $\mathcal{G}$ *permissible* if it is indexed by a set $Y$ that is an analytic subset of a compact metric space, such that $\mathcal{G} = \{g(\cdot, y), y \in Y\}$, and $g(\cdot, \cdot): \mathbb{R}^n \times Y \to \mathbb{R}$ is $\mathcal{L}_n \otimes \mathcal{B}(Y)$-measurable, where $\mathcal{L}_n$ is the class of Lebesgue measurable sets in $\mathbb{R}^n$ and $\mathcal{B}(Y)$ is the class of Borel subsets of $Y$. The metric on $Y$ will be important only insofar as it determines $\mathcal{B}(Y)$.

Let $(S, d)$ be a set $S$ equipped with a pseudometric $d$ and let $\varepsilon > 0$. A set $U \subset S$ is called an $\varepsilon$-*net* if each point in $S$ is within a $d$-distance at most $\varepsilon$ of some point in $U$.

We can now define *metric entropy*, a valuable concept introduced by Kolmogorov. Metric entropy is often also called $\varepsilon$-entropy, but we need to reserve the letter $\varepsilon$ for a different purpose. Accordingly, we define the $t$-*covering number* $N(t, S, d)$ of $(S, d)$ as the least cardinality of all $t$-nets. In other words, $N(t, S, d)$ is the least number of balls of radius $t$ with respect to $d$ that cover $S$. Then $H(t, S, d) = \log N(t, S, d)$ is called the $t$-*entropy* of $(S, d)$, and we can drop the argument $d$ when there is no possibility of confusion. This notion will mainly be used for subsets of $\mathcal{G}$. For $k \in \mathbb{N}$, we define a pseudonorm $|\cdot|_k$ on $\mathcal{G}$ by

$$|g|_k = \left(\frac{1}{k} \sum_{i=1}^{k} g(x_i)^2\right)^{1/2}, \qquad g \in \mathcal{G}.$$

Note that this pseudonorm depends on $x_1, \ldots, x_k$. For $\varepsilon > 0$, let

$$\mathcal{G}_k(\varepsilon, g_0) = \{g \in \mathcal{G} : |g - g_0|_k \leq \varepsilon\}.$$

Then we denote by $H(t, \mathcal{G}_k(\varepsilon, g_0))$ the $t$-entropy of $\mathcal{G}_k(\varepsilon, g_0)$ with respect to the pseudometric generated by the pseudonorm $|\cdot|_k$; again, this depends on $x_1, \ldots, x_k$. If $\mathcal{G}$ is a cone (i.e., $\mathcal{G} = s\mathcal{G}$ for all $s > 0$), then

(12) $\qquad H(t, \mathcal{G}_k(\varepsilon, g_0)) = H(st, \mathcal{G}_k(s\varepsilon, sg_0)) = H(st, s\mathcal{G}_k(\varepsilon, g_0)),$

for any $s > 0$. This follows from the fact that the balls of radius $t$ (with respect to $|\cdot|_k$) with centers $g_1, \ldots, g_m$ form a minimal cover of $\mathcal{G}_k(\varepsilon, g_0)$ if and only if the balls of radius $st$ with centers $sg_1, \ldots, sg_m$ form a minimal cover of $\mathcal{G}_k(s\varepsilon, sg_0)$. A *local entropy integral* $J_k(\varepsilon, \mathcal{G})$ can be defined for $a > 0$ and $0 < \varepsilon < 2^6 a$ by

(13) $\qquad J_k(\varepsilon, \mathcal{G}) = \max\left\{\int_{\varepsilon^2/(2^6 a)}^{\varepsilon} H(t, \mathcal{G}_k(\varepsilon, g_0))^{1/2} dt, \varepsilon\right\}.$

Note that this integral depends on $g_0$, $a$ and $x_1, \ldots, x_k$, although this is not explicit in the notation.

To state the principal technical result, a little more terminology is needed. The random variables $X_i$ are called *uniformly sub-Gaussian* if there are constants $A$ and $\tau$ such that, for $i = 1, 2, \ldots$, we have

(14) $\qquad\qquad\qquad A^2(E[e^{|X_i|^2/A^2}] - 1) \leq \tau^2.$

Note that if $X_i$ is a normal $N(0, \sigma^2)$ random variable for $i = 1, 2, \ldots$, then this condition is satisfied when $A = \tau = 2\sigma$.

The following result is due to van de Geer [32], Theorem 9.1 (see also [31]).

PROPOSITION 4.1. *Let $a > 0$ and let $X_i$, $i = 1, 2, \ldots$, be uniformly sub-Gaussian independent random variables satisfying* (14), *each with mean zero. Let $\mathcal{G}$ be a permissible class of real-valued functions on a subset $E$ of $\mathbb{R}^n$, let $g_0 \in \mathcal{G}$, and let $(x_i)$ be a sequence in $E$. Let $J_k(\varepsilon, \mathcal{G})$ be defined by* (13), *and suppose that $\Psi$ is a function with $\Psi(\varepsilon) \geq J_k(\varepsilon, \mathcal{G})$ for all $k \in \mathbb{N}$ and such that $\Psi(\varepsilon)/\varepsilon^2$ is decreasing for $0 < \varepsilon < 2^6 a$. Then there is a constant $c = c(A, \tau)$ such that, for any $k \in \mathbb{N}$ and any $\varepsilon_k > 0$ with $\sqrt{k}\varepsilon_k^2 \geq c\Psi(\varepsilon_k)$, we have*

$$(15) \qquad \Pr(|g_0 - \hat{g}_k|_k > \varepsilon_k) \leq c e^{-k \varepsilon_k^2 / c^2} + \Pr\left(\frac{1}{k} \sum_{i=1}^{k} X_i^2 > a^2\right),$$

*for any least squares estimator $\hat{g}_k$ of $g_0$ with respect to $\mathcal{G}$ based on measurements at $x_1, \ldots, x_k$.*

It is crucial that the constant $c$ depends neither on $a$ nor on $k$. In Theorem 9.1 of [32], the fact that $c$ is independent of $k$ is not explicitly stated and requires some explanation. In our notation, the proof of Theorem 9.1 of [32] arrives at the inequality $\sqrt{k}\varepsilon_k \geq 16 C \Psi(\varepsilon_k)$, where $C$ is a constant independent of $k$. The assumptions *and* (13) yield

$$\sqrt{k} \geq 16 C \frac{\Psi(\varepsilon_k)}{\varepsilon_k^2} \geq 16 C \frac{J_k(\varepsilon_k, \mathcal{G})}{\varepsilon_k^2} \geq \frac{16 C}{\varepsilon_k},$$

or $\sqrt{k}\varepsilon_k \geq 16 C$. This allows the finite sum on the last line of page 149 of [32] to be bounded above by a geometric series whose sum depends on $k$ only in the required exponential form. (See the proof of Lemma 3.2 of [32] for a similar argument.)

The following result is implicit in pages 187–188 of [32]. A proof is provided for the convenience of the reader and because we need some details about the constants involved.

COROLLARY 4.2. *Suppose that the assumptions on the random variables $X_i$ and class $\mathcal{G}$ in Proposition* 4.1 *hold. For all $k \in \mathbb{N}$, let $\hat{g}_k$ be a least squares estimator of $g_0$ with respect to $\mathcal{G}$, based on measurements at $x_1, \ldots, x_k$. If there are positive constants $\alpha$, $t_0$ and $M$ such that*

$$(16) \qquad\qquad H(t, \mathcal{G}_k(\varepsilon_0, g_0)) \leq M^2 t^{-\alpha},$$

for all $k \in \mathbb{N}$, $0 < t \leq t_0$ and $\varepsilon_0 = 2^{13/2}\tau$, then, almost surely, there are constants $C = C(A, \tau, \alpha)$ and $N = N(A, \tau, \alpha, t_0, M)$ such that

$$(17) \quad |g_0 - \hat{g}_k|_k \leq \begin{cases} CM^{2/(2+\alpha)} k^{-1/(2+\alpha)}, & \text{if } \alpha < 2, \\ Ck^{-1/4} \log k, & \text{if } \alpha = 2, \\ CM^{1/\alpha} k^{-1/(2\alpha)}, & \text{if } \alpha > 2, \end{cases}$$

for $k \geq N$.

PROOF. Let $J_k(\varepsilon, \mathcal{G})$ be defined by (13) with $a = \sqrt{2}\tau$. We may suppose that $\tau > 0$ and therefore that $a > 0$. As $H(t, \mathcal{G}_k(\varepsilon, g_0))$ is an increasing function of $\varepsilon$ (with $t$ fixed), (16) holds when $\varepsilon_0$ is replaced by any $0 < \varepsilon \leq 2^6 a = \varepsilon_0$.

Consider first the case $\alpha < 2$ and let $0 < \varepsilon < 2^6 a$. For $0 < \varepsilon < t_0$, we have

$$\int_{\varepsilon^2/(2^6 a)}^{\varepsilon} H(t, \mathcal{G}_k(\varepsilon, g_0))^{1/2} \, dt \leq \frac{2M}{2-\alpha} \left( \varepsilon^{1-\alpha/2} - \left( \frac{\varepsilon^2}{2^6 a} \right)^{1-\alpha/2} \right).$$

As $H(t, \mathcal{G}_k(\varepsilon, g_0))$ is a decreasing function of $t$ (with $\varepsilon$ fixed), $\varepsilon \geq t_0$ implies

$$\int_{\varepsilon^2/(2^6 a)}^{\varepsilon} H(t, \mathcal{G}_k(\varepsilon, g_0))^{1/2} \, dt$$

$$= \int_{\varepsilon^2/(2^6 a)}^{t_0} H(t, \mathcal{G}_k(\varepsilon, g_0))^{1/2} \, dt + \int_{t_0}^{\varepsilon} H(t, \mathcal{G}_k(\varepsilon, g_0))^{1/2} \, dt$$

$$\leq \frac{2M}{2-\alpha} t_0^{1-\alpha/2} + H(t_0, \mathcal{G}_k(\varepsilon, g_0))^{1/2} (\varepsilon - t_0).$$

Let

$$(18) \quad \Psi(\varepsilon) = \begin{cases} \max\left\{ \dfrac{2M}{2-\alpha} \varepsilon^{1-\alpha/2}, \varepsilon \right\}, & \text{if } 0 < \varepsilon < t_0, \\ \max\left\{ \dfrac{2M}{2-\alpha} t_0^{1-\alpha/2} + M t_0^{-\alpha/2}(\varepsilon - t_0), \varepsilon \right\}, & \text{if } \varepsilon \geq t_0. \end{cases}$$

Then $\Psi(\varepsilon) \geq J(\varepsilon, \mathcal{G})$ and by (18), $\Psi(\varepsilon)/\varepsilon^2$ is a decreasing function of $\varepsilon > 0$. Suppose that $c > 0$ and let $\varepsilon_k = A_1 k^{-1/(2+\alpha)}$. If both

$$(19) \quad A_1 = \left( \frac{2Mc}{2-\alpha} \right)^{2/(2+\alpha)}$$

and

$$(20) \quad k > \max\left\{ \left( \frac{c}{A_1} \right)^{2(2+\alpha)/\alpha}, \left( \frac{A_1}{t_0} \right)^{2+\alpha} \right\}$$

hold, then one can check that $\varepsilon_k < t_0$ and (using this also) that $\sqrt{k}\varepsilon_k^2 \geq c\Psi(\varepsilon_k)$. As noted by van de Geer [32], page 150, (14) implies that

$$\Pr\left( \frac{1}{k} \sum_{i=1}^{k} X_i^2 > 2\tau^2 \right) \leq e^{-k\tau^2/(12A^2)}.$$

Thus, (15) yields

$$\Pr(|g_0 - \hat{g}_k|_k > A_1 k^{-1/(2+\alpha)}) \leq c e^{-A_1^2 k^{\alpha/(2+\alpha)}/c^2} + e^{-k\tau^2/(12A^2)}, \quad (21)$$

provided (19) and (20) hold. The sum over $k$ of the right-hand side of (21) converges, so by the Borel–Cantelli lemma we have, almost surely,

$$|g_0 - \hat{g}_k|_k \leq A_1 k^{-1/(2+\alpha)} = CM^{2/(2+\alpha)} k^{-1/(2+\alpha)},$$

say, for sufficiently large $k$. Therefore, (17) is true when $\alpha < 2$.

The argument when $\alpha \geq 2$ is similar; we omit the details. If $\alpha > 2$, we take

$$\Psi(\varepsilon) = \begin{cases} \max\left\{\dfrac{2M}{\alpha-2}\left(\dfrac{\varepsilon^2}{2^6 a}\right)^{1-\alpha/2}, \varepsilon\right\}, & \text{if } 0 < \varepsilon < t_0, \\ \max\left\{\dfrac{2M}{\alpha-2}\left(\dfrac{t_0^2}{2^6 a}\right)^{1-\alpha/2} + M t_0^{-\alpha/2}(\varepsilon - t_0), \varepsilon\right\}, & \text{if } \varepsilon \geq t_0, \end{cases}$$

and $\varepsilon_k = A_2 k^{-1/2\alpha}$ for a suitable constant $A_2$. If $\alpha = 2$, we can take

$$\Psi(\varepsilon) = \begin{cases} \max\{M(\log(2^6 a) - \log \varepsilon), \varepsilon\}, & \text{if } 0 < \varepsilon < t_0, \\ \max\{M(\log(2^6 a) - \log t_0 + (\varepsilon - t_0)/t_0), \varepsilon\}, & \text{if } \varepsilon \geq t_0, \end{cases}$$

and $\varepsilon_k = A_3 k^{-1/4} \log k$ for a suitable constant $A_3$. In both cases it can be checked that $\varepsilon_k < t_0$ and then that $\sqrt{k}\varepsilon_k^2 \geq c\Psi(\varepsilon_k)$ when $k$ is sufficiently large. (The case $\alpha = 2$ is qualitatively different, as $A_3$ can be chosen independent of $M$.) □

**5. Least squares estimation of support functions.** Suppose that $K$ is an unknown convex body in $\mathbb{R}^n$, and $(u_i)$ is a sequence in $S^{n-1}$. For $k \in \mathbb{N}$, the support function $h_K$ of $K$ is measured at $u_1, u_2, \ldots, u_k$. The measurements

$$y_i = h_K(u_i) + X_i, \quad (22)$$

$i = 1, 2, \ldots, k$, are noisy, the $X_i$'s being independent random variables with zero mean and finite variance. We want to find a convex body with the property that its support function values at $u_1, \ldots, u_k$ best approximate the measurements $y_1, \ldots, y_k$. In order to apply the results of the previous section, we let $E = S^{n-1}$ and

$$\mathcal{G} = \{h_L : L \in \mathcal{K}^n\},$$

the class of support functions, throughout this section.

LEMMA 5.1.  *The class $\mathcal{G}$ is permissible.*

PROOF. Referring to the above definition of the term "permissible," the index set $Y = \mathcal{K}^n$ is a subset of the family $\mathcal{F}_n$ of all closed subsets of $\mathbb{R}^n$. The latter, endowed with the hit-and-miss topology, is a compact metrizable space; see, for example, [29], Satz 1.1.1. By Satz 1.3.2 of [29], $Y$ is a Borel set in $\mathcal{F}_n$, so it is analytic in $\mathcal{F}_n$. Although the induced topology on $Y$ as a subset of $\mathcal{F}_n$ is coarser than the topology induced by the Hausdorff metric, the respective families of Borel sets coincide; see [29], Satz 1.3.2. The mapping $(K, u) \mapsto h_K(u)$ is continuous with respect to both topologies, so the parametrization mapping is Borel measurable. $\square$

Fix $k \in \mathbb{N}$ and $K \in \mathcal{K}^n$. In accordance with the notation of the previous section [see (11)], we let $(\widehat{h_K})_k$ be a least squares estimator for $h_K$ with respect to $\mathcal{G}$ based on measurements at $u_1, \ldots, u_k$, so that $h_K$ now plays the role of the function $g_0$. As $\mathcal{G}$ is a closed cone in the usual Banach space of continuous functions on the sphere [and the objective function in (11) is continuous on this space], a least squares estimator always exists. For $h \colon S^{n-1} \to \mathbb{R}$, the pseudonorm $|h|_k$ is now given by

$$(23) \qquad |h|_k = \left( \frac{1}{k} \sum_{i=1}^{k} h(u_i)^2 \right)^{1/2}.$$

The following lemma provides an upper bound for the $L_2$ distance between two convex bodies $L$ and $M$ contained in a ball $SB$ in terms of the pseudometric $|h_L - h_M|_k$.

LEMMA 5.2. *Let $S > 0$ and let $L$ and $M$ be convex bodies in $\mathbb{R}^n$ contained in $SB$. Let $\{u_1, \ldots, u_k\}$ be a subset of $S^{n-1}$. Then*

$$(24) \qquad \delta_2(L, M) \leq (k\omega_k)^{1/2}(|h_L - h_M|_k + 2\Delta_k S),$$

*where $\Delta_k$ and $\omega_k$ are defined by* (5) *and* (7), *respectively.*

PROOF. As in Section 3, denote the Voronoi cells corresponding to $\{u_1, \ldots, u_k\}$ by $\Omega_i$, $1 \leq i \leq k$. If $u \in \Omega_i$, we have $\|u - u_i\| \leq \Delta_k$ and hence,

$$h_L(u) \leq h_L(u_i) + h_L(u - u_i)$$
$$\leq h_L(u_i) + \|u - u_i\| h_L\left( \frac{u - u_i}{\|u - u_i\|} \right) \leq h_L(u_i) + \Delta_k S.$$

Similarly,

$$h_M(u_i) \leq h_M(u) + h_M(u_i - u) \leq h_M(u) + \Delta_k S.$$

Therefore,

$$h_L(u) - h_M(u) \leq h_L(u_i) - h_M(u_i) + 2\Delta_k S,$$

and interchanging $L$ and $M$, we obtain

$$|h_L(u) - h_M(u)| \leq |h_L(u_i) - h_M(u_i)| + 2\Delta_k S.$$

Therefore,

$$\delta_2(L,M)^2 = \int_{S^{n-1}} (h_L(u) - h_M(u))^2 \, du$$

$$\leq \sum_{i=1}^{k} \int_{\Omega_i} (|h_L(u_i) - h_M(u_i)| + 2\Delta_k S)^2 \, du$$

$$\leq \omega_k \sum_{i=1}^{k} (|h_L(u_i) - h_M(u_i)| + 2\Delta_k S)^2$$

$$\leq \omega_k \left( \left( \sum_{i=1}^{k} (h_L(u_i) - h_M(u_i))^2 \right)^{1/2} + \left( \sum_{i=1}^{k} (2\Delta_k S)^2 \right)^{1/2} \right)^2$$

$$= k\omega_k (|h_L - h_M|_k + 2\Delta_k S)^2. \qquad \square$$

We shall also need the next lemma, which under the assumption $K \subset RB$ and a mild condition on the sequence $(u_i)$ yields the radius of a ball containing $L$ in terms of the pseudometric $|h_K - h_L|_k$.

LEMMA 5.3. *Let $K$ and $L$ be convex bodies in $\mathbb{R}^n$. Suppose that $K \subset RB$ for some $R > 0$, and that $(u_i)$ is an evenly spread sequence in $S^{n-1}$. Then there are constants $C_0 = C_0((u_i)) > 0$ and $N_0 = N_0((u_i)) \in \mathbb{N}$ such that*

$$L \subset (C_0|h_K - h_L|_k + 2R)B,$$

*for all $k \geq N_0$.*

PROOF. Fix $k$ and choose $x_k \in L$, where we may assume that $\|x_k\| > 2R$ since otherwise $L \subset 2RB$. Then $h_L(u) \geq x_k \cdot u$ for all $u \in S^{n-1}$. Let $v_k = x_k/\|x_k\|$. Choose $t_0 > 0$ small enough that, for each $u \in S^{n-1}$ and any $v, w \in C_{t_0}(u)$, we have $v \cdot w \geq 1/2$. (Of course, $t_0$ does not depend on $k$.) If $u_i \in C_{t_0}(v_k)$, then

$$|h_K(u_i) - h_L(u_i)| \geq x_k \cdot u_i - R > \frac{\|x_k\|}{2} - R > 0.$$

Therefore,

$$|h_K - h_L|_k^2 \geq \frac{1}{k} \sum_{u_i \in C_{t_0}(v_k)} |h_K(u_i) - h_L(u_i)|^2$$

$$\geq \left( \frac{\|x_k\|}{2} - R \right)^2 \frac{1}{k} |\{u_1, \ldots, u_k\} \cap C_{t_0}(v_k)|$$

CONVERGENCE OF RECONSTRUCTION ALGORITHMS 17$$\geq c\left(\frac{\|x_k\|}{2} - R\right)^2,$$

for some $c > 0$ and all $k \geq N$, say, because $(u_i)$ is evenly spread. [Note that $c$ and $N$ depend only on $(u_i)$.] From this, we obtain

$$\|x_k\| \leq 2\left(\frac{1}{\sqrt{c}}|h_K - h_L|_k + R\right),$$

for $k \geq N$, and the result follows. $\square$

The following result is due to Bronshtein [3]. His definition of entropy uses $\log_2$ instead of the natural logarithm, requiring an extra constant factor in the lower bound.

PROPOSITION 5.4. *Let $\mathcal{K}^n(B)$ denote the space of compact convex subsets of the unit ball $B$ in $\mathbb{R}^n$, endowed with the Hausdorff metric. Then for $0 < t < 10^{-12}/(n-1)$, the $t$-entropy $H(t, \mathcal{K}^n(B))$ of $\mathcal{K}^n(B)$ satisfies*

(25) $$\frac{\kappa_{n-1} \log 2}{8^{n-1}(n-1)} t^{-(n-1)/2} \leq H(t, \mathcal{K}^n(B)) \leq (\log 12) 10^6 n^{5/2} t^{-(n-1)/2}.$$

Let $\varepsilon$ and $t$ be positive numbers and let $k \in \mathbb{N}$. As before, let

$$\mathcal{G}_k(\varepsilon, h_K) = \{h_L \in \mathcal{G} : |h_K - h_L|_k \leq \varepsilon\},$$

and let

$$H(t, \mathcal{G}_k(\varepsilon, h_K)) = H(t, \mathcal{G}_k(\varepsilon, h_K), |\cdot|_k)$$

be the $t$-entropy of $\mathcal{G}_k(\varepsilon, h_K)$ with respect to the pseudometric generated by $|\cdot|_k$.

COROLLARY 5.5. *Let $(u_i)$ be an evenly spread sequence in $S^{n-1}$ and let $K$ be a convex body in $\mathbb{R}^n$ with $K \subset RB$ for some $R > 0$. Then there are constants $t_1 = t_1(n, (u_i))$ and $C_1 = C_1(n, (u_i))$ such that*

(26) $$H(t, \mathcal{G}_k(\varepsilon, h_K)) \leq C_1 R^{(n-1)/2} t^{-(n-1)/2},$$

*for all $k \in \mathbb{N}$, $0 < \varepsilon \leq R$ and $0 < t \leq Rt_1$.*

PROOF. We first make the following claim: There is a constant $s_0 = s_0(n, (u_i)) > 0$ such that, for all $k \in \mathbb{N}$ and $h_L \in \mathcal{G}_k(\varepsilon, h_K)$, there is an $L' \in (R/s_0)\mathcal{K}^n(B)$ with $h_{L'}(u_i) = h_L(u_i)$, for $i = 1, \ldots, k$.

The claim will be proved later. Assuming it is true, we observe that if $h_L \in (s_0/R)\mathcal{G}_k(\varepsilon, h_K)$, then there is an $L' \in \mathcal{K}^n(B)$ such that $|h_L - h_M|_k =$ $|h_{L'} - h_M|_k$ for any compact convex set $M$ in $\mathbb{R}^n$. It follows from this and (12) with $s = s_0/R$ that

$$H(t, \mathcal{G}_k(\varepsilon, h_K)) = H(s_0 t/R, (s_0/R)\mathcal{G}_k(\varepsilon, h_K))$$
$$\leq H(s_0 t/R, \{h_{L'} : L' \in \mathcal{K}^n(B)\}, |\cdot|_k).$$

Since

$$|h_L - h_M|_k \leq \|h_L - h_M\|_\infty = \delta(L, M),$$

for any two compact convex subsets $L$ and $M$ in $\mathbb{R}^n$, we have

$$H(t, \mathcal{G}_k(\varepsilon, h_K)) \leq H(s_0 t/R, \{h_{L'} : L' \in \mathcal{K}^n(B)\}, |\cdot|_k) \leq H(s_0 t/R, \mathcal{K}^n(B)).$$

Now (26) is an immediate consequence of (25) if we put $t_1 = 2 \cdot 10^{-12}/((n-1)s_0)$ and $C_1 = (\log 12)10^6 n^{5/2}(2s_0)^{(n-1)/2}$.

It remains to prove the claim. Let $h_L \in \mathcal{G}_k(\varepsilon, h_K)$. By Lemma 5.3, there are constants $C_0 = C_0((u_i))$ and $N_0 = N_0((u_i))$ such that if $k \geq N_0$, then

$$(27) \qquad L \subset (C_0 \varepsilon + 2R)B \subset (C_0 + 2)RB.$$

For such $k$, we let $L' = L$. Now let $k \leq N_0$. Since $h_L \in \mathcal{G}_k(\varepsilon, h_K)$, we have

$$(28) \qquad |h_L(u_i)|^2 \leq k|h_L|_k^2 \leq N_0(\varepsilon + |h_K|_k)^2 \leq 4N_0 R^2,$$

for $i = 1, \ldots, k$. Let $I \subset \{1, \ldots, N_0\}$ be nonempty, and consider the continuous function $f_I$ on $S^{n-1}$ defined by

$$f_I(u) = \sum_{i \in I} |u \cdot u_i|.$$

For $u$ in the span of $\{u_i : i \in I\}$, $f_I(u) > 0$. Therefore, we can choose $a_0 = a_0(n, (u_i)) > 0$ such that, for any such $I$ and $u$ in the span of $\{u_i : i \in I\}$, $f_I(u) \geq a_0$.

Suppose that $\{u_1, \ldots, u_k\}$ spans $\mathbb{R}^n$. The polyhedron $P = \bigcap_{i=1}^k \{x \in \mathbb{R}^n : x \cdot u_i \leq h_L(u_i)\}$ satisfies $h_P(u_i) = h_L(u_i)$ for $i = 1, \ldots, k$, but may be unbounded. Let $L' = \text{conv}\{x_1, \ldots, x_m\}$, where $x_1, \ldots, x_m$ are the vertices of $P$. Then $L'$ is bounded and satisfies $h_{L'}(u_i) = h_L(u_i)$ for $i = 1, \ldots, k$. Any vertex $x_j$ of $P$ is an intersection of $n$ hyperplanes with linearly independent normals $u_{i_1}, \ldots, u_{i_n}$, say. Using (28) with $L$ replaced by $L'$, we obtain, for any $x_j \neq o$,

$$\|x_j\| a_0 \leq \|x_j\| \sum_{p=1}^n \left| \frac{x_j}{\|x_j\|} \cdot u_{i_p} \right| = \sum_{p=1}^n |x_j \cdot u_{i_p}| \leq \sum_{p=1}^n h_{L'}(u_{i_p}) \leq 2n\sqrt{N_0} R.$$

Thus, $L' \subset (2n\sqrt{N_0}/a_0)RB$.

If the span of $\{u_1, \ldots, u_k\}$ is a proper subspace $S$ of $\mathbb{R}^n$, the above argument can be applied to $L|S$ regarded as a compact convex subset of $S$ to obtain the same inclusion. In view of (27), which holds for all $k \geq N_0$ with $L$ replaced by $L'$, we conclude that $L' \subset (R/s_0)B$ for all $k \in \mathbb{N}$, where $s_0 = \min\{1/(C_0 + 2), a_0/(2n\sqrt{N_0})\}$ depends only on $n$ and $(u_i)$. □

CONVERGENCE OF RECONSTRUCTION ALGORITHMS 19THEOREM 5.6. *Let $(u_i)$ be an evenly spread sequence in $S^{n-1}$ and let $X_i$, $i = 1, 2, \ldots$, be uniformly sub-Gaussian independent random variables satisfying* (14), *each with mean zero. Let $K$ be a convex body in $\mathbb{R}^n$ with $K \subset RB$, where $R \geq 2^{13/2}\tau$, and for $k \in \mathbb{N}$, let $\widehat{(h_K)}_k$ be a least squares estimator of $h_K$ with respect to $\mathcal{G}$, based on measurements at $u_1, \ldots, u_k$. Then, almost surely, there are constants $C_2 = C_2(A, \tau, n)$ and $N_2 = N_2(A, \tau, n, R, (u_i))$ such that*

$$(29) \quad |h_K - \widehat{(h_K)}_k|_k \leq \begin{cases} C_2 R^{(n-1)/(n+3)} k^{-2/(n+3)}, & \text{if } n = 2, 3, 4, \\ C_2 k^{-1/4} \log k, & \text{if } n = 5, \\ C_2 R^{1/2} k^{-1/(n-1)}, & \text{if } n \geq 6, \end{cases}$$

*for $k \geq N_2$.*

PROOF. Let $\varepsilon_0 = 2^{13/2}\tau$. Since $0 < \varepsilon_0 \leq R$, Corollary 5.5 yields

$$(30) \qquad H(t, \mathcal{G}_k(\varepsilon_0, h_K)) \leq C_1 R^{(n-1)/2} t^{-(n-1)/2},$$

for $k \in \mathbb{N}$ and $0 < t \leq Rt_1$. By (30), we may apply Corollary 4.2, with $\alpha = (n-1)/2$, $t_0 = Rt_1$ and

$$M^2 = C_1 R^{(n-1)/2},$$

to conclude that (17) holds, almost surely, with $M$ as above and $C = C(A, \tau, \alpha) = C_2(A, \tau, n)$ and $N = N(A, \tau, \alpha, t_0, M) = N_2(A, \tau, n, R, (u_i))$. □

COROLLARY 5.7. *Let $(u_i)$ be an evenly spread sequence in $S^{n-1}$ and let $X_i$, $i = 1, 2, \ldots$, be independent $N(0, \sigma^2)$ random variables. Let $K$ be a convex body in $\mathbb{R}^n$ with $K \subset RB$, where $R \geq \sigma^{15/2}$, and for $k \in \mathbb{N}$, let $\widehat{(h_K)}_k$ be a least squares estimator of $h_K$ with respect to $\mathcal{G}$, based on measurements at $u_1, \ldots, u_k$. Then, almost surely, there are constants $C_3 = C_3(n)$ and $N_3 = N_3(\sigma, n, R, (u_i))$ such that*

$$|h_K - \widehat{(h_K)}_k|_k \leq \begin{cases} C_3 \sigma^{4/(n+3)} R^{(n-1)/(n+3)} k^{-2/(n+3)}, & \text{if } n = 2, 3, 4, \\ C_3 \sigma k^{-1/4} \log k, & \text{if } n = 5, \\ C_3 \sigma^{1/2} R^{1/2} k^{-1/(n-1)}, & \text{if } n \geq 6, \end{cases}$$

(31)
*for $k \geq N_3$.*

PROOF. As was noted earlier, we may take $A = \tau = 2\sigma$ in Theorem 5.6 and conclude that if $K \subset RB$ and $R \geq \sigma^{15/2}$, then, almost surely, the least squares estimators $\widehat{(h_K)}_k$ for $h_K$ with respect to $\mathcal{G}$ satisfy (29), where the dependence of $C_2$ and $N_2$ on $A$ and $\tau$ is replaced by dependence on $\sigma$. Instead we now use scaled measurements

$$\lambda y_i = \lambda h_K(u_i) + \lambda X_i,$$

with some $\lambda > 0$, to estimate the support function $h_{\lambda K} = \lambda h_K$ of the scaled convex body $\lambda K$. Then $\lambda(\widehat{h_K})_k$ is a least squares estimator for $h_{\lambda K}$. Also, $\lambda X_i$, $i = 1, 2, \ldots$, are independent normal $N(0, (\lambda\sigma)^2)$ random variables. Replacing $K$, $R$ and $\sigma$ by $\lambda K$, $\lambda R$ and $\lambda\sigma$, respectively, we conclude that, almost surely, there are constants $c_0 = c_0(\lambda\sigma, n)$ and $n_0 = n_0(\lambda\sigma, n, \lambda R, (u_i))$ such that

$$\tag{32} |\lambda h_K - \lambda(\widehat{h_K})_k|_k \leq c_0(\lambda R)^{b_n} f_n(k),$$

for $k \geq n_0$, where $R^{b_n}$ and $f_n(k)$ are the functions of $R$ and $k$, respectively, in (29). When $\lambda = 1/\sigma$, (32) becomes

$$\tag{33} |h_K - (\widehat{h_K})_k|_k \leq C_3 \sigma^{1-b_n} R^{b_n} f_n(k),$$

where $C_3 = C_3(n)$ and $k \geq N_3 = N_3(\sigma, n, R, (u_i))$. Substituting $b_n$ and $f_n$ from (29) into (33), we arrive at (31). $\square$

**6. Convergence of the Prince–Willsky algorithm.** Let $u_1, \ldots, u_k$ be fixed vectors in $S^{n-1}$ whose positive hull is $\mathbb{R}^n$. We say that the nonnegative real numbers $h_1, \ldots, h_k$ are *consistent* if there is a compact convex set $L$ in $\mathbb{R}^n$ such that $h_L(u_i) = h_i$, $i = 1, \ldots, k$. If $h_1, \ldots, h_k$ are consistent, there will be many such sets $L$; let $P(h_1, \ldots, h_k)$ denote the one that is the polytope defined by

$$\tag{34} P(h_1, \ldots, h_k) = \bigcap_{i=1}^{k} \{x \in \mathbb{R}^n : x \cdot u_i \leq h_i\}.$$

For $n = 2$ and vectors $u_1, \ldots, u_k$ equally spaced in $S^1$, the following algorithm was proposed and implemented by Prince and Willsky [27].

ALGORITHM NOISYSUPPORTLSQ.

*Input*: Natural numbers $n \geq 2$ and $k \geq n + 1$; vectors $u_i \in S^{n-1}$, $i = 1, \ldots, k$, whose positive hull is $\mathbb{R}^n$; noisy support function measurements

$$y_i = h_K(u_i) + X_i,$$

$i = 1, \ldots, k$, of an unknown convex body $K$ in $\mathbb{R}^n$, where the $X_i$'s are independent $N(0, \sigma^2)$ random variables.

*Task*: Construct a convex polytope $\hat{P}_k$ in $\mathbb{R}^n$ that approximates $K$, with facet outer normals belonging to the set $\{u_i : i = 1, \ldots, k\}$.

*Action*: Solve the following constrained linear least squares problem (LLS1):

$$\tag{35} \min_{h_1, \ldots, h_k} \sum_{i=1}^{k} (y_i - h_i)^2,$$

$$\tag{36} \text{subject to } h_1, \ldots, h_k \text{ are consistent}.$$

Let $\hat{h}_1, \ldots, \hat{h}_k$ be a solution of (LLS1) and let $\hat{P}_k = P(\hat{h}_1, \ldots, \hat{h}_k)$. $\square$

Naturally any implementation of Algorithm NoisySupportLSQ involves making explicit the constraint (36). Although we do not need to address this problem for our purposes, a few remarks are appropriate. For $n = 2$, this was done by Prince and Willsky [27] for vectors $u_1, \ldots, u_k$ equally spaced in $S^1$, and by Lele, Kulkarni and Willsky [21] for arbitrary vectors $u_1, \ldots, u_k$, by means of an inequality constraint of the form $Ah \leq 0$, where $h = (h_1, \ldots, h_k)$ and $A$ is a certain matrix. For general $n$, this is more difficult and was studied by Karl, Kulkarni, Verghese and Willsky [19]. (In these papers there is no mention of Rademacher's condition for consistency when $n = 2$, or of Firey's extension (see [28], page 47) of Rademacher's condition to $n \geq 2$.) The authors of [19] did not implement the algorithm for $n \geq 3$; an implementation for $n = 3$ and certain special sets of directions was carried out by Gregor and Rannou [14].

If the positive hull of $\{u_1, \ldots, u_k\}$ is not $\mathbb{R}^n$, then (34) could still be considered as output of the Algorithm NoisySupportLSQ, if consistency of $h_1, \ldots, h_k$ is extended to closed convex sets which may be unbounded. We choose not to do this, however. Indeed, if $(u_i)$ is a dense sequence of vectors in $S^{n-1}$, then, for sufficiently large $k$, the positive hull of $u_1, \ldots, u_k$ is $\mathbb{R}^n$ and in this case, Algorithm NoisySupportLSQ produces a polytope $\hat{P}_k$ as output. We now establish conditions under which $\hat{P}_k$ converges, almost surely, to $K$ as $k \to \infty$. Of course, the denseness of $(u_i)$ is a necessary condition for such convergence.

The following theorem establishes the strong consistency of Algorithm NoisySupportLSQ when $(u_i)$ is evenly spread.

THEOREM 6.1. *Let $K$ be a convex body in $\mathbb{R}^n$ and let $(u_i)$ be an evenly spread sequence in $S^{n-1}$. If $\hat{P}_k$ is an output from Algorithm NoisySupportLSQ as stated above, then, almost surely,*

$$\lim_{k \to \infty} \delta(K, \hat{P}_k) = 0.$$

PROOF. Theorem 5.6 and $(\widehat{h_K})_k = h_{\hat{P}_k}$ imply that, almost surely, we have

(37) $$\lim_{k \to \infty} |h_K - h_{\hat{P}_k}|_k = 0.$$

Fix a realization for which (37) holds.

By Lemma 5.3, there is an $S > 0$ such that $\hat{P}_k \subset SB$ for all $k$. According to Blaschke's selection theorem, the set $\{\hat{P}_1, \hat{P}_2, \ldots\}$ is relatively compact in the space of convex bodies in $\mathbb{R}^n$ with the Hausdorff metric. To prove $\lim_{k \to \infty} \hat{P}_k = K$, it is therefore enough to show that $K$ is the only accumulation point of $(\hat{P}_k)$.

Let $\tilde{K}$ be an arbitrary accumulation point of this sequence. Then a subsequence of $(h_{\hat{P}_k})$ converges uniformly to $h_{\tilde{K}}$. This and (37) can be applied to the right-hand side of

$$|h_K - h_{\tilde{K}}|_k \le |h_K - h_{\hat{P}_k}|_k + |h_{\hat{P}_k} - h_{\tilde{K}}|_k$$

to show that a subsequence $(|h_K - h_{\tilde{K}}|_{k'})$ converges to 0. For each $k'$ in this subsequence,

$$|h_K - h_{\tilde{K}}|_{k'} = \|h_K - h_{\tilde{K}}\|_{L_2(\mu_{k'})}$$

is the $L_2$ norm of $h_K - h_{\tilde{K}}$ with respect to the probability measure $\mu_{k'}$ that assigns a mass $1/k'$ to each of the points $u_1, \ldots, u_{k'}$. As the set of probability measures in $S^{n-1}$ is weakly compact, there is a subsequence $(\mu_{k''})$ of $(\mu_{k'})$ that converges weakly to a probability measure $\mu$. Using the continuity of support functions, we conclude that

$$0 = \lim_{k \to \infty} \|h_K - h_{\tilde{K}}\|_{L_2(\mu_{k''})} = \|h_K - h_{\tilde{K}}\|_{L_2(\mu)}.$$

We claim that, since $(u_i)$ is evenly spread, the support of $\mu$ is $S^{n-1}$; this will then imply $h_K = h_{\tilde{K}}$ and, hence, $K = \tilde{K}$. To prove the claim, suppose that $G$ is a nonempty open set in $S^{n-1}$ such that $\mu(G) = 0$. Choose an open cap $C_t(u) \subset G$, $t > 0$, and a nonnegative continuous real-valued function $f$ on $S^{n-1}$ with support contained in $G$ and such that $f \ge 1$ on $C_t(u)$. Then the fact that $(u_i)$ is evenly spread implies that

$$0 < \liminf_{k \to \infty} \int_{S^{n-1}} \mathbb{1}_{C_t(u)}(v) \, d\mu_k(v) \le \lim_{k \to \infty} \int_{S^{n-1}} f(v) \, d\mu_k(v)$$
$$= \int_{S^{n-1}} f(v) \, d\mu(v) \le \|f\|_\infty \mu(G) = 0,$$

where $\mathbb{1}_{C_t(u)}$ denotes the characteristic function of $C_t(u)$. This contradiction completes the proof. $\square$

The conclusion of the following theorem is stronger than that of Theorem 6.1 since it provides convergence rates. However, the hypothesis on the sequence $(u_i)$ is also stronger; see Lemma 3.2, which also guarantees the existence of suitable sequences $(u_i)$.

THEOREM 6.2. *Let $\sigma > 0$ and let $K$ be a convex body in $\mathbb{R}^n$ such that $K \subset RB$ for some $R \ge 2^{15/2}\sigma$. Let $(u_i)$ be a sequence in $S^{n-1}$ with $\Delta_k = O(k^{-1/(n-1)})$. If $\hat{P}_k$ is an output from Algorithm NoisySupportLSQ as stated above, then, almost surely, there are constants $C_4 = C_4(n, (u_i))$ and $N_4 = N_4(\sigma, n, R, (u_i))$ such that*

$$(38) \quad \delta_2(K, \hat{P}_k) \le \begin{cases} C_4 \sigma^{4/(n+3)} R^{(n-1)/(n+3)} k^{-2/(n+3)}, & \text{if } n = 2, 3, 4, \\ C_4 \sigma k^{-1/4} \log k, & \text{if } n = 5, \\ C_4 (R + (\sigma R)^{1/2}) k^{-1/(n-1)}, & \text{if } n \ge 6, \end{cases}$$

for $k \geq N_4$.

Also, there are constants $C_5 = C_5(n, (u_i))$ and $N_5 = N_5(\sigma, n, R, (u_i))$ such that

$$(39) \quad \delta(K, \hat{P}_k) \leq \begin{cases} C_5 \sigma^{8/((n+1)(n+3))} R^{(n-1)(n+5)/((n+1)(n+3))} \\ \quad \times k^{-4/((n+1)(n+3))}, & \text{if } n = 2, 3, 4, \\ C_5 \sigma^{1/3} R^{2/3} k^{-1/12} (\log k)^{1/3}, & \text{if } n = 5, \\ C_5 (R + \sigma^{1/(n+1)} R^{n/(n+1)}) k^{-2/(n^2-1)}, & \text{if } n \geq 6, \end{cases}$$

for $k \geq N_5$.

PROOF. By Lemma 5.3, we can apply Lemma 5.2 with $L = K$, $M = \hat{P}_k$ and $S = C_0 |h_K - h_{\hat{P}_k}|_k + 2R$ to obtain

$$(40) \quad \delta_2(K, \hat{P}_k) \leq (k\omega_k)^{1/2}(4\Delta_k R + (2\Delta_k C_0 + 1)|h_K - h_{\hat{P}_k}|_k),$$

for all $k \geq N_0$. By Lemma 3.2, $k\omega_k = O(1)$ and we also have $\widehat{(h_K)}_k = h_{\hat{P}_k}$. The various estimates for $\delta_2(K, \hat{P}_k)$ now follow from the corresponding estimates for $|h_K - \widehat{(h_K)}_k|_k$ in Corollary 5.7.

To obtain the estimates for $\delta(K, \hat{P}_k)$, we combine those just found and the relation (1) that yields

$$\delta(K, \hat{P}_k) \leq cS^{(n-1)/(n+1)} \delta_2(K, \hat{P}_k)^{2/(n+1)},$$

where $S = C_0 |h_K - h_{\hat{P}_k}|_k + 2R$, for all $k \geq N_0$. $\square$

**7. Reconstruction from brightness function measurements.** Suppose that $K$ is an unknown origin-symmetric convex body in $\mathbb{R}^n$, and $(u_i)$ is a sequence of directions in $S^{n-1}$. For $k \in \mathbb{N}$, the brightness function $b_K$ of $K$ is measured at $u_1, u_2, \ldots, u_k$. The measurements

$$(41) \quad y_i = b_K(u_i) + X_i,$$

$i = 1, 2, \ldots, k$, are noisy, the $X_i$'s being independent random variables with zero mean and finite variance. We want to find an origin-symmetric convex body with the property that its brightness function values at $u_1, \ldots, u_k$ best approximate the measurements $y_1, \ldots, y_k$.

The following algorithm was proposed by Gardner and Milanfar [12]. Since it is convenient for us to describe it in somewhat different language, we briefly explain how it works in the case of exact measurements, a situation analyzed in detail by Gardner and Milanfar [13]. The algorithm proceeds in two phases, motivated by the connection between zonoids, projection bodies and surface area measures outlined in Section 2. In the first phase, a constrained least squares problem is solved to find a zonotope $Z$ with $h_Z(u_i) = b_K(u_i)$, $i = 1, \ldots, k$. This zonotope is the projection body of a polytope whose surface area measure can easily be calculated from $Z$. The second phase reconstructs the polytope from this known surface area measure.

ALGORITHM NOISYBRIGHTLSQ.

*Input*: Natural numbers $n \geq 2$ and $k$; vectors $u_i \in S^{n-1}$, $i = 1, \ldots, k$, that span $\mathbb{R}^n$; noisy brightness function measurements

$$y_i = b_K(u_i) + X_i,$$

$i = 1, \ldots, k$, of an unknown origin-symmetric convex body $K$ in $\mathbb{R}^n$, where the $X_i$'s are independent normal $N(0, \sigma^2)$ random variables.

*Task*: Construct a convex polytope $\hat{Q}_k$ in $\mathbb{R}^n$ that approximates $K$.

*Action*:

*Phase* I: Find a zonotope $\hat{Z}_k \in \mathcal{Z}^n$ that solves the following least squares problem:

$$\min_{Z \in \mathcal{Z}^n} \sum_{i=1}^{k} (y_i - h_Z(u_i))^2. \tag{42}$$

Calculate the (finitely supported) surface area measure $S(\hat{Q}_k, \cdot)$ of the origin-symmetric polytope $\hat{Q}_k$ satisfying

$$\hat{Z}_k = \Pi \hat{Q}_k. \tag{43}$$

*Phase* II: Reconstruct $\hat{Q}_k$ from $S(\hat{Q}_k, \cdot)$ (or directly from $\hat{Z}_k$, if possible). □

It was observed by Gardner and Milanfar [13] that the remark after Proposition 2.1 can be used in Phase I; this shows that a zonotope $\hat{Z}_k$ solving (42) exists. Moreover, as $\hat{Z}_k$ can be assumed to be a sum of line segments, each parallel to a node corresponding to $U = \{u_1, \ldots, u_k\}$, only the length of these line segments has to be determined. [Note, however, that this restriction on the direction of the line segments is not required in (42).] This turns (42) into a finite-dimensional quadratic program which can be solved using standard software. When $n = 2$, Phase II is simple (see [13], page 284, but note also the statistically improved method proposed by Poonawala, Milanfar and Gardner [26]). For $n \geq 3$, Phase II is highly nontrivial, but can be performed by means of the previously known algorithm MinkData (see [13] for references).

When the brightness function measurements are exact, it was proved by Gardner and Milanfar ([13], Theorem 6.1) that if $(u_i)$ is dense in $S^{n-1}$, then the outputs $\hat{Q}_k$ [corresponding to the first $k$ directions in $(u_i)$] converge to $K$, as $k \to \infty$. For a convergence proof that applies to noisy measurements, we can apply our results from Section 5. We begin with a suitable form of Lemma 5.3. Recall definition (9) of the symmetrized sequence $(u_i^*)$.

LEMMA 7.1. *Let $K$ and $L$ be origin-symmetric convex bodies in $\mathbb{R}^n$. Suppose that $K \subset RB$ for some $R > 0$, and that $(u_i)$ is a sequence in $S^{n-1}$* such that $(u_i^*)$ is evenly spread. Then there are constants $C_0^* = C_0^*((u_i)) > 0$ and $N_0^* = N_0^*((u_i)) \in \mathbb{N}$ such that

$$L \subset (C_0^*|h_K - h_L|_k + 2R)B,$$

for all $k \geq N_0^*$.

PROOF. This follows easily from Lemma 5.3, if $(u_i)$ and $k$ are replaced by $(u_i^*)$ and $2k$, respectively, and

$$\frac{1}{2k}\sum_{i=1}^{2k}(h_K(u_i^*) - h_L(u_i^*))^2 = |h_K - h_L|_k^2$$

is taken into account. $\square$

The next theorem gives the strong consistency of Algorithm NoisyBrightLSQ when $(u_i^*)$ is evenly spread.

THEOREM 7.2. *Let $K$ be a convex body in $\mathbb{R}^n$ and let $(u_i)$ be a sequence in $S^{n-1}$ such that $(u_i^*)$ is evenly spread. If $\hat{Q}_k$ is an output from Algorithm NoisyBrightLSQ as stated above, then, almost surely,*

(44) $$\lim_{k \to \infty} \delta(K, \hat{Q}_k) = 0.$$

PROOF. Choose $0 < r < R$ such that $rB \subset K \subset RB$. Then $\Pi(rB) \subset \Pi K \subset \Pi(RB)$, so

(45) $$sB \subset \Pi K \subset tB,$$

where $s = \kappa_{n-1}r^{n-1}$ and $t = \kappa_{n-1}R^{n-1}$. Theorem 5.6 and $(\widehat{h_{\Pi K}})_k = h_{\hat{Z}_k}$ imply that, almost surely, we have

(46) $$\lim_{k \to \infty} |h_{\Pi K} - h_{\hat{Z}_k}|_k = 0.$$

Fix a realization for which (46) holds.

By (45), (46) and Lemma 7.1 with $K$, $L$ and $R$ replaced by $\Pi K$, $\hat{Z}_k$ and $\kappa_{n-1}R^{n-1}$, respectively, there is an $S > 0$ such that $\hat{Z}_k \subset SB$ holds for all $k$. We can now apply Blaschke's selection theorem and the argument used in the proof of Theorem 6.1 to conclude that

(47) $$\lim_{k \to \infty} \delta(\Pi K, \hat{Z}_k) = 0,$$

as $(u_i^*)$ is evenly spread in $S^{n-1}$.

When $n = 2$, it is easy to see that $\Pi K$ and $\Pi \hat{Q}_k$ are rotations about the origin by $\pi/2$ of $2K$ and $2\hat{Q}_k$, respectively. (See, e.g., [10], Theorem 4.1.4.) Therefore, (44) follows immediately from (47).

Suppose that $n \geq 3$. By (45) and (47), we have

$$\text{(48)} \qquad \frac{s}{2}B \subset \hat{Z}_k = \Pi \hat{Q}_k \subset \frac{3}{2}tB,$$

for sufficiently large $k$. (Note that the fact that $\hat{Z}_k$ is $n$-dimensional for sufficiently large $k$ guarantees the existence of $\hat{Q}_k$.) Exactly the same argument as in the proof of Lemma 4.2 of [13] [beginning with formula (16) in that paper] leads from (48) to

$$\text{(49)} \qquad r_0 B \subset \hat{Q}_k \subset R_0 B,$$

for sufficiently large $k$, where

$$\text{(50)} \qquad R_0 = \frac{3n\kappa_n}{\kappa_{n-1}}\left(\frac{3}{2}\right)^{1/(n-1)} \frac{R^n}{r^{n-1}} \quad \text{and} \quad r_0 = \frac{\kappa_{n-1}r^{n-1}}{2^n R_0^{n-2}}.$$

Since $rB \subset K \subset RB$ implies $r_0 B \subset K \subset R_0 B$, we can apply (47) and Proposition 2.2 with $L = \hat{Q}_k$ to obtain (44). □

The results from Section 5 also give rates of convergence. However, we are able to do better, at least for $3 \leq n \leq 5$, by replacing the class of regression functions by the smaller family

$$\tilde{\mathcal{G}} = \{h_Z : Z \in \mathcal{Z}^n\}.$$

Note that this class is permissible, since it is easy to check that $\mathcal{Z}^n$ is a Borel set in $\mathcal{K}^n$.

In the plane, the class $\mathcal{Z}^2(B)$ of origin-symmetric zonoids contained in $B$ is just the class of origin-symmetric convex bodies contained in $B$. Using this fact, an appropriate modification of the proof of [3], Theorem 4, of the lower bound in (25) can be made that shows there is a constant $c > 0$ such that

$$\text{(51)} \qquad H(t, \mathcal{Z}^2(B)) \geq ct^{-1/2},$$

for sufficiently small $t > 0$. It follows that the exact entropy exponents for $\mathcal{Z}^2(B)$ and $\mathcal{K}^2(B)$ are the same, namely, $-1/2$. For $n \geq 3$, however, the following theorem represents a dramatic improvement.

THEOREM 7.3. *Let $\mathcal{Z}^n(B)$ denote the space of origin-symmetric zonoids contained in the unit ball $B$ in $\mathbb{R}^n$, endowed with the Hausdorff metric. If $n \geq 3$, then for all $0 < t < 1/2$ and any $\eta > 0$,*

$$\text{(52)} \qquad H(t, \mathcal{Z}^n(B)) = O(t^{-2(n-1)/(n+2)-\eta}).$$

PROOF. Let $t > 0$. Suppose that $K$ is a zonoid in $\mathcal{Z}^n(B)$. Clearly, there are an $N = N(n, t) \in \mathbb{N}$, depending only on $n$ and $t$, and a zonotope $Z$ such that

$$(53) \qquad K \subset Z \subset (1 + t/2)K,$$

where $Z = \sum_{i=1}^{N} a[-v_i, v_i]$, for some $0 < a < 1$ and $v_i \in S^{n-1}$, $i = 1, \ldots, N$.

Let $S$ be a $t/(4N)$-net in $[0, 1]$ and let $U$ be a $t/(4N)$-net in $S^{n-1}$. Let $s$ be the closest point in $S$ to $a$, let $u_i$ be the closest point in $U$ to $v_i$, $i = 1, \ldots, N$, and let

$$(54) \qquad Z' = \sum_{i=1}^{N} s[-u_i, u_i].$$

For each $i = 1, \ldots, N$, $L_i = a[-v_i, v_i]$ and $M_i = s[-u_i, u_i]$ are origin-symmetric line segments whose Hausdorff distance apart is bounded by the distance between the points $av_i$ and $su_i$. Using this, we obtain

$$\delta(Z, Z') = \|h_Z - h_{Z'}\|_\infty \le \sum_{i=1}^{N} \|h_{L_i} - h_{M_i}\|_\infty \le \sum_{i=1}^{N} \|av_i - su_i\|$$

$$\le \sum_{i=1}^{N} (\|av_i - au_i\| + \|au_i - su_i\|) \le \sum_{i=1}^{N} \left(\frac{t}{4N} + \frac{t}{4N}\right) = \frac{t}{2}.$$

From this and (53), we obtain $\delta(K, Z') \le t$.

By Proposition 3.1, we can choose $S$ and $U$ so that $|S| = O(N/t)$ and $|U| = O((N/t)^{n-1})$. With this choice, the number of zonotopes of the form (54) is $O((N/t)^{nN})$. Therefore, the $t$-entropy of $\mathcal{Z}^n(B)$ is

$$(55) \qquad H(t, \mathcal{Z}^n(B)) = O\left(N \log\left(\frac{N}{t}\right)\right).$$

Bourgain and Lindenstrauss [2] proved that, for $0 < t < 1/2$, one can take

$$(56) \qquad N = N(3, t) = O\left(t^{-4/5}\left(\log \frac{1}{t}\right)^{2/5}\right)$$

when $n = 3$ and

$$(57) \qquad N = N(4, t) = O\left(t^{-1}\left(\log \frac{1}{t}\right)^{3/2}\right)$$

when $n = 4$. They also obtained a good bound for $n \ge 5$, but this was improved by Matoušek [24], who obtained

$$(58) \qquad N = N(n, t) = O(t^{-2(n-1)/(n+2)}),$$

for $n \ge 5$. Substituting (56), (57) and (58) into (55), we obtain (52). □

Let $\varepsilon$ and $t$ be positive numbers and let $k \in \mathbb{N}$ and $Z \in \mathcal{Z}^n$ be given. In accordance with earlier notation, let

$$\tilde{\mathcal{G}}_k(\varepsilon, h_Z) = \{h_L \in \tilde{\mathcal{G}} : |h_Z - h_L|_k \leq \varepsilon\},$$

and let $H(t, \tilde{\mathcal{G}}_k(\varepsilon, h_Z)) = H(t, \tilde{\mathcal{G}}_k(\varepsilon, h_Z), |\cdot|_k)$ be the $t$-entropy of $\tilde{\mathcal{G}}_k(\varepsilon, h_Z)$ with respect to the pseudometric generated by $|\cdot|_k$.

COROLLARY 7.4. *Let $(u_i)$ be a sequence in $S^{n-1}$ such that $(u_i^*)$ is evenly spread and let $Z$ be an origin-symmetric zonoid in $\mathbb{R}^n$ with $Z \subset RB$ for some $R > 0$. If $n \geq 3$, then for any $\eta > 0$, there are constants $t_6 = t_6(n, (u_i), \eta)$ and $C_6 = C_6(n, (u_i), \eta)$ such that*

$$(59) \qquad H(t, \tilde{\mathcal{G}}_k(\varepsilon, h_Z)) \leq C_6 R^{2(n-1)/(n+2)+\eta} t^{-2(n-1)/(n+2)-\eta},$$

*for all $k \in \mathbb{N}$, $0 < \varepsilon \leq R$ and $0 < t \leq Rt_6$.*

PROOF.  This follows from Theorem 7.3, exactly as Corollary 5.5 follows from Proposition 5.4.  □

Lemma 3.3 guarantees the existence of sequences satisfying the hypothesis of the next theorem.

THEOREM 7.5. *Let $\sigma > 0$ and let $K$ be an origin-symmetric convex body in $\mathbb{R}^n$ such that $K \subset RB$, where $\kappa_{n-1} R^{n-1} \geq 2^{15/2} \sigma$. Let $(u_i)$ be a sequence of directions in $S^{n-1}$ with $\Delta_k^* = O(k^{-1/(n-1)})$, and suppose that $\hat{Z}_k$ is a corresponding solution of (42). If $n = 2$, then, almost surely, there are constants $C_7 = C_7((u_i))$ and $N_7 = N_7(\sigma, R, (u_i))$ such that*

$$(60) \qquad \delta_2(\Pi K, \hat{Z}_k) \leq C_7 \sigma^{4/5} R^{1/5} k^{-2/5},$$

*for $k \geq N_7$. If $n = 3$ or $4$, there is a constant $\gamma_0 = \gamma_0(n) > 0$ such that if $0 < \gamma < \gamma_0$, then, almost surely, there are constants $C_8 = C_8(n, (u_i), \gamma)$ and $N_8 = N_8(\sigma, n, R, (u_i), \gamma)$ such that*

$$(61) \quad \delta_2(\Pi K, \hat{Z}_k) \leq C_8 \sigma^{(n+2)/(2n+1)-\gamma} R^{(n-1)^2/(2n+1)+\gamma} k^{-(n+2)/(4n+2)+\gamma},$$

*for $k \geq N_8$.*

*Finally, if $n \geq 5$, there are constants $C_9 = C_9(n, (u_i))$ and $N_9 = N_9(\sigma, n, R, (u_i))$ such that*

$$(62) \qquad \delta_2(\Pi K, \hat{Z}_k) \leq C_9 R^{n-1} k^{-1/(n-1)},$$

*for $k \geq N_9$.*

PROOF. Let $\varepsilon_0 = 2^{15/2}\sigma$ and $\eta > 0$. As $K \subset RB$, we have $\Pi K \subset \kappa_{n-1} R^{n-1} B$. Since $0 < \varepsilon_0 \leq \kappa_{n-1} R^{n-1}$, Corollary 5.5 (for $n = 2$, using $\tilde{\mathcal{G}} \subset \mathcal{G}$) and Corollary 7.4 (for $n \geq 3$) with $Z$ and $R$ replaced by $\Pi K$ and $\kappa_{n-1} R^{n-1}$, respectively, yield

$$H(t, \tilde{\mathcal{G}}_k(\varepsilon, h_{\Pi K})) \leq \tilde{C}_6 R^{(n-1)\alpha} t^{-\alpha}, \tag{63}$$

for all $k \in \mathbb{N}$ and $0 < t \leq Rt_6$, where

$$\alpha = \begin{cases} 1/2, & \text{if } n = 2, \\ 2(n-1)/(n+2) + \eta, & \text{if } n \geq 3. \end{cases} \tag{64}$$

If $\eta < 6/(n+2)$, then $\alpha < 2$, so applying Corollary 4.2 with this $\alpha$, $t_0 = Rt_6$, $M^2 = \tilde{C}_6 R^{(n-1)\alpha}$ and $\mathcal{G}$ replaced by $\tilde{\mathcal{G}}$, we conclude that, almost surely, there are constants $c_0 = c_0(\sigma, n, \eta)$ and $n_0 = n_0(\sigma, n, R, (u_i), \eta)$ such that

$$|h_{\Pi K} - h_{\hat{Z}_k}|_k \leq c_0 R^{(n-1)\alpha/(2+\alpha)} k^{-1/(2+\alpha)}, \tag{65}$$

for $k \geq n_0$.

The dependence on $\sigma$ is dealt with by the device used in proving Corollary 5.7. By using scaled measurements,

$$\lambda y_i = \lambda b_K(u_i) + \lambda X_i = \lambda h_{\Pi K}(u_i) + \lambda X_i = h_{\Pi(\lambda^{1/(n-1)} K)}(u_i) + \lambda X_i,$$

replacing $K$, $R$ and $\sigma$ by $\lambda^{1/(n-1)} K$, $\lambda^{1/(n-1)} R$ and $\lambda\sigma$, respectively, and then setting $\lambda = 1/\sigma$, we obtain from (65) the inequality

$$|h_{\Pi K} - h_{\hat{Z}_k}|_k \leq c_1 \sigma^{2/(2+\alpha)} R^{(n-1)\alpha/(2+\alpha)} k^{-1/(2+\alpha)}, \tag{66}$$

which holds, almost surely, for some $c_1 = c_1(n, \eta)$ and $k \geq n_0$.

By Lemma 3.3, we may apply Lemma 7.1, with $K$, $L$ and $R$ replaced by $\Pi K$, $\hat{Z}_k$ and $\kappa_{n-1} R^{n-1}$, respectively, to obtain $\Pi K, \hat{Z}_k \subset SB$ for all $k \geq N_0^*$, where

$$S = C_0^* |h_{\Pi K} - h_{\hat{Z}_k}|_k + 2\kappa_{n-1} R^{n-1}. \tag{67}$$

In Lemma 5.2 we make similar substitutions and replace the set $\{u_1, \ldots, u_k\}$ by $\{u_1^*, \ldots, u_{2k}^*\}$, to conclude that

$$\delta_2(\Pi K, \hat{Z}_k) \leq (2k\omega_k^*)^{1/2}(|h_{\Pi K} - h_{\hat{Z}_k}|_k + 2\Delta_k^* S),$$

for all $k \geq N_0^*$. This, (66) and the fact that by Lemma 3.3 we have $k\omega_k^* = O(1)$ imply that there are constants $C' = C'(n, (u_i), \eta) > 0$ and $N' = N'(\sigma, n, R, (u_i), \eta) > 0$ such that

$$\delta_2(\Pi K, \hat{Z}_k) \leq C'(\sigma^{2/(2+\alpha)} R^{(n-1)\alpha/(2+\alpha)} k^{-1/(2+\alpha)} + R^{n-1} k^{-1/(n-1)}) \tag{68}$$

for all $k \geq N'$. For $n \geq 5$ and large $k$, the second term dominates and (62) follows. For $n \leq 4$ and large $k$, the first term dominates; then (64) and (68) yield (60) for $n = 2$ and (61) for $n = 3$ and $4$. □

The next theorem gives rates of convergence for Algorithm NoisyBrightLSQ in terms of the Hausdorff metric. For $n \geq 3$, we omit the dependence on $R$ because this is complicated by the use of Proposition 2.2; as we mentioned above, no particular effort was made to obtain optimal results in the estimate (4).

THEOREM 7.6. *Let $\sigma > 0$ and let $K$ be a convex body in $\mathbb{R}^n$ such that $K \subset RB$, where $\kappa_{n-1}R^{n-1} \geq 2^{15/2}\sigma$. Let $(u_i)$ be a sequence of directions in $S^{n-1}$ with $\Delta_k^* = O(k^{-1/(n-1)})$, and suppose that $\hat{Q}_k$ is an output of Algorithm NoisyBrightLSQ as stated above. If $n = 2$, then, almost surely, there are constants $C_{10} = C_{10}(\sigma, (u_i))$ and $N_{10} = N_{10}(\sigma, R, (u_i))$ such that*

$$\delta(K, \hat{Q}_k) \leq C_{10} R^{7/15} k^{-4/15}, \tag{69}$$

*for $k \geq N_{10}$.*

*If $n \geq 3$, suppose, in addition, that $rB \subset K$ for some $0 < r < R$. For $n = 3$ or $4$, there is a constant $\gamma_1 = \gamma_1(n) > 0$ such that if $0 < \gamma < \gamma_1$, then, almost surely, there are constants $C_{11} = C_{11}(\sigma, n, r, R, (u_i), \gamma)$ and $N_{11} = N_{11}(\sigma, n, r, R, (u_i), \gamma)$ such that*

$$\delta(K, \hat{Q}_k) \leq C_{11} k^{-(n+2)/(n(n+4)(2n+1))+\gamma}, \tag{70}$$

*for $k \geq N_{11}$.*

*If $n \geq 5$ and $\gamma > 0$, then, almost surely, there are constants $C_{12} = C_{12}(\sigma, n, r, R, (u_i), \gamma)$ and $N_{12} = N_{12}(\sigma, n, r, R, (u_i), \gamma)$ such that*

$$\delta(K, \hat{Q}_k) \leq C_{12} k^{-2/((n-1)n(n+4))+\gamma}, \tag{71}$$

*for $k \geq N_{12}$.*

PROOF. Suppose that $n = 2$. Then $\Pi K$ and $\Pi \hat{Q}_k$ are rotations about the origin by $\pi/2$ of $2K$ and $2\hat{Q}_k$, respectively. Then (69) follows directly from (60) and (1).

Now suppose that $n \geq 3$. We have $sB \subset \Pi K \subset tB$, where $s = \kappa_{n-1} r^{n-1}$ and $t = \kappa_{n-1} R^{n-1}$. Note that (61) (for $n = 3$ or $4$), (62) (for $n \geq 5$) and (1) imply that, almost surely, there is a constant $N_{13} = N_{13}(\sigma, n, r, R, (u_i))$ such that (48) holds for all $k \geq N_{13}$. As in the proof of Theorem 7.2, we can conclude that

$$r_0 B \subset K, \hat{Q}_k \subset R_0 B, \tag{72}$$

for $k \geq N_{13}$, where $r_0$ and $R_0$ are given by (50). The desired results, (70) for $n = 3$ or $4$ and (71) for $n \geq 5$, now follow from Proposition 2.2 (with $L = \hat{Q}_k$) and Theorem 7.5. $\square$

The use of Proposition 2.2 in the previous theorem introduces a factor that worsens the convergence rates considerably. For example, when $n = 3$, we obtain a convergence rate of approximately $k^{-1/30}$!

**8. Monte Carlo simulations.** The theory of empirical processes that underlies our theoretical results suggests that the rates of convergence obtained in Corollary 5.7, for support function estimation with respect to the pseudonorm $|\cdot|_k$, are suboptimal for $n \geq 5$ (cf. page 162 of [32]). However, for $n \leq 4$, we expect them to be optimal, and this should carry over to the (identical) rates for Algorithm NoisySupportLSQ with respect to the $L_2$ metric, given in Theorem 6.2, as well as to the rates obtained in connection with Algorithm BrightLSQ, given in Theorem 7.5. On the other hand, we cannot expect the rates given in Theorems 6.2 and 7.6 involving the Hausdorff metric to be optimal, in view of the use of (1) (and, in the case of Theorem 7.6, the use of Proposition 2.2).

Extensive Monte Carlo simulations were run. The simulations are restricted to the case $n = 2$, since there does not appear to be a fully satisfactory implementation of Algorithm NoisySupportLSQ in higher dimensions (see the remarks in Section 6) and our present implementation of Algorithm NoisyBrightLSQ is too slow to allow enough iterations (we hope to improve this in the near future). In each simulation, a polygon was reconstructed 1000 times from noisy measurements of its support function or brightness function, using our implementations of Algorithm NoisySupportLSQ or Algorithm NoisyBrightLSQ, respectively. We developed the computer programs with the help of Chris Eastman, Greg Richardson, Thomas Riehle and Chris Street (work done as Western Washington University undergraduates) and Amyn Poonawala (at UC Santa Cruz).

Before describing the results of the simulations, we need to clarify the role of $R$ and the assumption in the above theorems about its relation to the noise level $\sigma$. For example, the inequality $R \geq 2^{15/2}\sigma$ is often assumed in order to prove that $d(K, \hat{P}_k) \leq C\sigma^a R^b f_n(k)$, where $d$ is the pseudometric $|\cdot|_k$, the $L_2$ metric or the Hausdorff metric. To test the dependence on $k$ or on $\sigma$ over any fixed range $\sigma_0 \leq \sigma \leq \sigma_1$, we can obviously choose $R$ large enough so that $K \subset RB$ and $R \geq 2^{15/2}\sigma_1$ is satisfied. We claim that the condition $R \geq 2^{15/2}\sigma$ also does not play any essential role in testing the dependence on $R$, and that we can view $R$ as a scaling factor of $K$. To see this, suppose $K$, $\sigma$ and a range $0 < \lambda_0 \leq \lambda \leq \lambda_1$ of scaling factors are given. Choose $R_0$ large enough so that $K \subset R_0 B$ and $\lambda_0 R_0 \geq 2^{15/2}\sigma$. Then $\lambda K \subset (\lambda R_0)B$ and $\lambda R_0 \geq 2^{15/2}\sigma$ for $\lambda \geq \lambda_0$. Replacing $K$ and $R$ in our theorems by $\lambda K$ and $\lambda_0 R_0$, respectively, we obtain

$$d(\lambda K, \widehat{P(\lambda)}_k) \leq C\sigma^a(\lambda R_0)^b f_n(k) = C\sigma^a R_0^b \lambda^b f_n(k) = C'\lambda^b f_n(k),$$

where $\widehat{P(\lambda)}_k$ is the output polytope for input $\lambda K$ and where $C'$ does not depend on $\lambda$. Thus, the exponent for $\lambda$ is the same as that for $R$ above, proving the claim.

Two input polygons were used, the regular 11-gon and irregular 9-gon displayed in Figure 1. Some results for the regular 11-gon are shown in Figure 2.

Each graph shows the results from 1000 iterations of Algorithm NoisySupportLSQ. The graphs are divided vertically into two groups of six graphs, corresponding to noise levels 0.1 and 1. In the left-hand column, the error (i.e., the distance between the input polygon and output polygon) is measured with the pseudonorm $|\cdot|_k$, while in the middle and right-hand columns, the $L_2$ and Hausdorff distances, respectively, are used instead. Each graph shows a curve giving the average error over all 1000 iterations, and points plotted above the curve giving the maximum error over the 1000 iterations. In each group of six graphs the top row shows error against the scaling factor $R$ varying from $R=0.2$ to $R=6$ in steps of 0.2, where the support function is always measured in the 35 directions at angles $0, 2\pi/35, 4\pi/35, \ldots, 68\pi/35$. The second row in each group of six graphs shows error against the number $k$ of measurements [in directions at angles $0, 2\pi/k, 4\pi/k, \ldots, 2(k-1)\pi/k$] varying from 20 to 100 in steps of 5, with the scaling factor $R$ fixed at 1.

For each of the 12 graphs in Figure 2, we used standard software to fit a curve of the form $CR^b$ or $Ck^c$ (for error against $R$ or $k$, resp.) to the points representing the averages over the 1000 iterations, and we repeated this for the points representing the maxima over the 1000 iterations. The corresponding values of $b$ and $c$ are shown in Table 1.

The case $n=2$ of Corollary 5.7 and Theorem 6.2 suggests that the appropriate values are $b = 1/5 = 0.2$ and $c = -2/5 = -0.4$ when errors are measured with $|\cdot|_k$ and the $L_2$ metric, and $b = 7/15 = 0.4666\ldots$ and $c = -4/15 = -0.2666\ldots$ when errors are measured with the Hausdorff metric. Of course, these theorems apply only for sufficiently large values of $k$ depending on both the noise level $\sigma$ and the scale factor $R$.

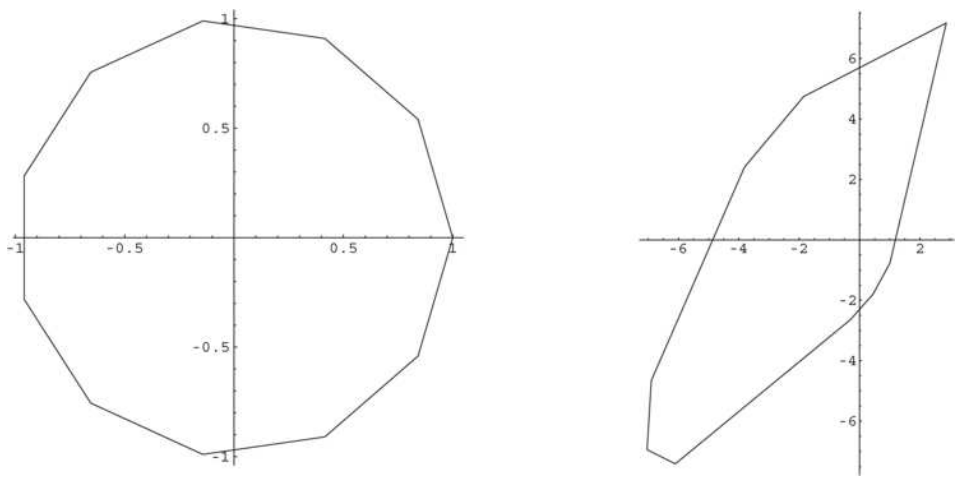

FIG. 1. *A regular* 11-*gon and irregular* 9-*gon.*

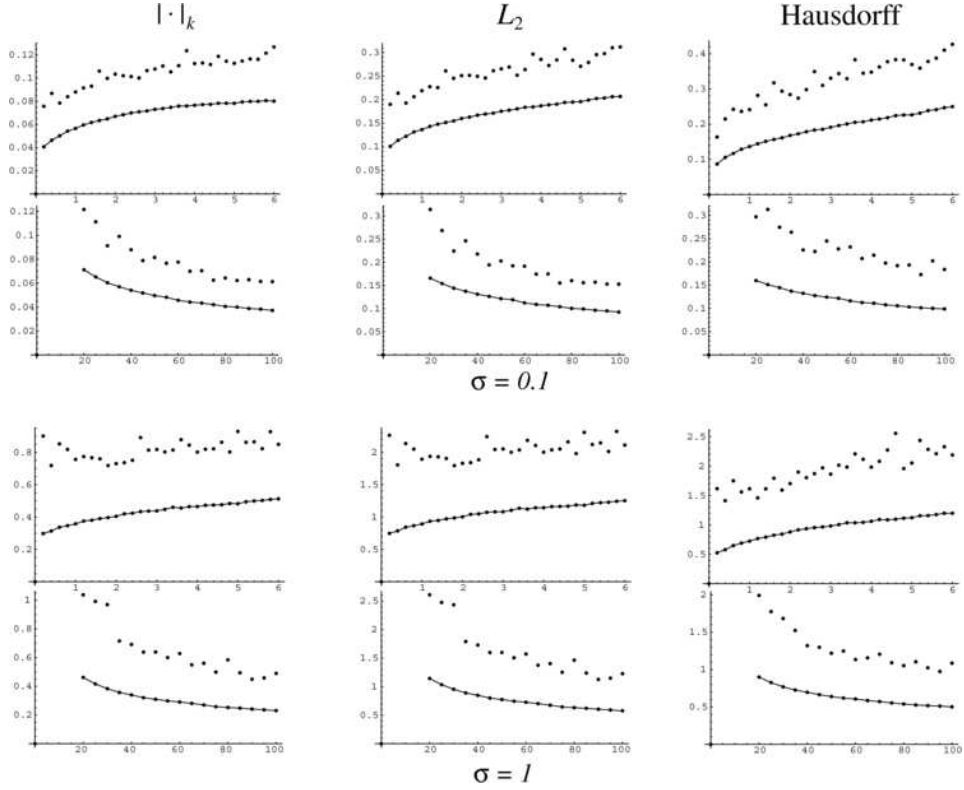

FIG. 2. *Error against $R$ and $k$ for the regular 11-gon.*

TABLE 1
*Fit for average and maximum error against $R$ and $k$ (11-gon)*

|  | | Average | | | Maximum | | |
|---|---|---|---|---|---|---|---|
|  | Error | $\|\cdot\|_k$ | $L_2$ | Hausdorff | $\|\cdot\|_k$ | $L_2$ | Hausdorff |
| $\sigma = 0.1$ | $b$ | 0.2020 | 0.2226 | 0.3248 | 0.1521 | 0.1567 | 0.2521 |
|  | $c$ | $-0.4006$ | $-0.3593$ | $-0.3052$ | $-0.4415$ | $-0.4468$ | $-0.3262$ |
| $\sigma = 1$ | $b$ | 0.1787 | 0.1684 | 0.2668 | 0.2771 | 0.2295 | 0.1686 |
|  | $c$ | $-0.4268$ | $-0.4202$ | $-0.3628$ | $-0.5338$ | $-0.5347$ | $-0.4316$ |

Despite the varying values in Table 1, we believe that the results of our Monte Carlo simulations are compatible with the expectations outlined in the first paragraph of this section, except perhaps in the case of Hausdorff error against scale. When the noise level is $\sigma = 0.1$, the values given in Table 1 for the $|\cdot|_k$ error, $b = 0.2020$ and $c = -0.4006$, for the average of the 1000 iterations are in very close agreement with theory, and the agreement is only slightly worse for the other metrics and at the high noise level $\sigma = 1$, except for Hausdorff error against scale.

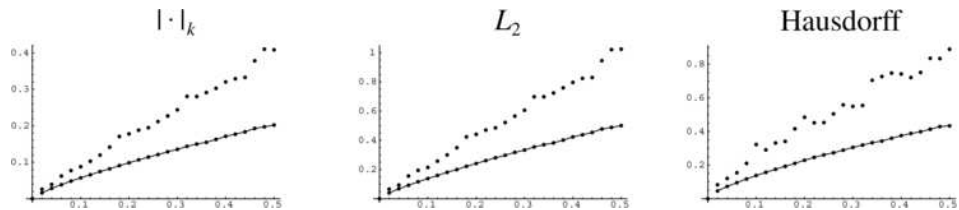

FIG. 3. *Error against $\sigma$ for the regular 11-gon.*

Naturally, the results for the maximum of the 1000 iterations are more unreliable due to the stochastic nature of the simulations. However, a poor fit does not necessarily contradict Corollary 5.7 and Theorem 6.2 (note the words "almost surely" in the statements of these theorems), especially for high noise levels. Better fits can be expected for data representing 1 or 2 standard deviations, for example, above the average for the 1000 iterations. For example, when $\sigma = 1$, the data representing 2 standard deviations above the average gives $c = -0.4505$, $-0.4479$ and $-0.3745$ for the $|\cdot|_k$, $L_2$ and Hausdorff errors, respectively (compare the three numbers at the right of the bottom row in Table 1).

In Figure 3 the three types of errors for the regular 11-gon are plotted against noise level $\sigma$ varying from $\sigma = 0.02$ to $\sigma = 0.5$ in steps of 0.02. Here the support function is always measured in the 35 directions at angles $0, 2\pi/35, 4\pi/35, \ldots, 68\pi/35$, and the polygon is unscaled (i.e., $R = 1$). As before, each graph shows a curve giving the average error over all 1000 iterations, and points plotted above the curve giving the maximum error over the 1000 iterations. The exponents $a$ for curves of best fit of the form $C\sigma^a$ are, for the average, $a = 0.7894$, $0.8038$ and $0.7150$ for the $|\cdot|_k$, $L_2$ and Hausdorff errors, respectively. The corresponding exponents for the maximum error are $a = 0.9286$, $0.9346$ and $0.7593$. These are in good agreement with the values $a = 4/5 = 0.8$ for the $|\cdot|_k$ and $L_2$ errors given in Corollary 5.7 and Theorem 6.2. The less convincing agreement with the value $a = 8/15 = 0.5333\ldots$ for the Hausdorff error given in Theorem 6.2 is not surprising, since the discrepancy between $L_2$ and Hausdorff errors that occurs via (1) is smaller for the regular 11-gon than for a general polygon.

Simulations for the regular 11-gon at a low noise level, $\sigma = 0.01$, as well as for the irregular 9-gon in Figure 1, were also compatible with theory. For the details, see [11], Section 8.

Suppose that we attempt to reconstruct an origin-symmetric planar convex body $K$, first with Algorithm NoisyBrightLSQ, using $k$ noisy brightness function measurements at angles $0, \pi/k, \ldots, (k-1)\pi/k$, and then with Algorithm NoisySupportLSQ, using $2k$ noisy support function measurements at angles $0, \pi/k, \ldots, (2k-1)\pi/k$. The two output polygons will, in general, be different, but apart from the noise, this is only because the two sets of measurements do not "match." Indeed, for any angle $\alpha \in [0, 2\pi)$,

$$h_K(\alpha \pm \pi/2) = b_K(\alpha)/2, \tag{73}$$

in view of the origin symmetry of $K$. In fact, there is a very close relationship between our implementations of Algorithms NoisySupportLSQ and NoisyBrightLSQ when $n = 2$. If we run Algorithm NoisyBrightLSQ with noisy brightness function values $y_i$ measured at angles $\alpha_i$, $i = 1, \ldots, k$, in the interval $[0, \pi)$, our implementation will produce an origin-symmetric output polygon $\hat{Q}_k$ with outer normals among the directions $\alpha_i \pm \pi/2$, $i = 1, \ldots, k$; see [13]. Using this fact and (73), it is easy to prove that if we then run Algorithm NoisySupportLSQ using $y_i/2$ as noisy support function value at angle $\alpha_i \pm \pi/2$, $i = 1, \ldots, k$, the output polygon will also be $\hat{Q}_k$. Thus, very similar results can be expected from the two algorithms when $n = 2$ and $K$ is origin symmetric, and we verified this by performing simulations of 1000 iterations of Algorithm NoisyBrightLSQ for a regular origin-symmetric 12-gon and an affinely regular origin-symmetric octagon. We omit the details, noting only that values of $a$ and $b$ indicated by the case $n = 2$ of Theorems 7.5 and 7.6 are the same as those above, and that the observed agreement was similar in all respects to that detailed above for Algorithm NoisySupportLSQ.

**9. Application to a stereological problem.** In this section the convergence results above are used to obtain strong consistency of an estimator for the directional measure of a random collection of fibers. Details about the following notions can be found in Chapter 9 of [30]. A *fiber* is a $C^1$ curve of finite length, and a *fiber process* $Y$ is a random element with values in the family of locally finite collections of fibers in $\mathbb{R}^n$. We assume that $Y$ is *stationary* (the term *homogeneous* is also used), meaning that the distribution of $Y$ is translation invariant. Suppose that $A$ is a Borel set in $\mathbb{R}^n$ with $V_n(A) = 1$ and $E$ is an origin-symmetric Borel set in $S^{n-1}$. Let $\mu(E)$ be the mean total length of the union of all fiber points in $A$ with a unit tangent vector in $E$. Due to the stationarity of $Y$, $\mu(E)$ is independent of $A$ and so this definition gives rise to a unique even Borel measure $\mu$ in $S^{n-1}$ called the *directional measure* of $Y$. We also assume that, almost surely, the fibers of $Y$ do not all lie in parallel hyperplanes, so that $\mu$ is not concentrated on a great sphere. The *length density* $\overline{L} = \mu(S^{n-1})$ is the mean total length of fibers per unit volume. The probability measure $\mu/\overline{L}$, called the *rose of directions*, can be interpreted as the distribution of a unit tangent vector at a "typical" fiber point, and hence, can be used to quantify anisotropy of $Y$.

In applications, the fiber process $Y$ often cannot be observed directly, but only via its intersections with planes. Due to the stationarity, we can restrict our considerations to hyperplanes containing the origin. For each $u \in S^{n-1}$, let $\gamma(u)$ be the mean number of points in $Y \cap u^\perp$ per unit $(n-1)$-dimensional volume. The function $\gamma$ is called the *rose of intersections* of $Y$. It is well known that

$$\gamma(u) = \int_{S^{n-1}} |u \cdot v| \, d\mu(v) \tag{74}$$

for all $u \in S^{n-1}$. As $h(u) = |u \cdot v|$, $u \in S^{n-1}$ is the support function of the line segment $[-v, v]$, (74) shows that $\gamma$ is the support function of a zonoid $Z$, called the *associated zonoid* or *Steiner compact* of $Y$. Minkowski's existence theorem implies that there is a convex body $K$ with surface area measure $2\mu$. As

$$\tfrac{1}{2} \int_{S^{n-1}} |u \cdot v| \, dS(K, v) = h_{\Pi K}(u), \tag{75}$$

for all $u \in S^{n-1}$ [see, e.g., [28], equation (5.3.34)], we have $h_Z(u) = \gamma(u) = h_{\Pi K}(u)$, $u \in S^{n-1}$, and so $Z = \Pi K$.

Since $\gamma(u) = h_{\Pi K}(u) = b_K(u)$ for $u \in S^{n-1}$ and $\mu = (1/2)S(K, \cdot)$, the following slightly modified version of Phase I of Algorithm NoisyBrightLSQ allows the reconstruction of an approximation $\hat{\mu}_k$ to $\mu$ from noisy measurements of $\gamma$.

ALGORITHM NOISYROSELSQ.

*Input*: Natural numbers $n \geq 2$ and $k$; vectors $u_i \in S^{n-1}$, $i = 1, \ldots, k$, that span $\mathbb{R}^n$; noisy measurements

$$y_i = \gamma(u_i) + X_i, \tag{76}$$

$i = 1, \ldots, k$, of the rose of intersections $\gamma$ of an unknown stationary fiber process $Y$ in $\mathbb{R}^n$, where the $X_i$'s are independent $N(0, \sigma^2)$ random variables.

*Task*: Construct a finitely supported measure $\hat{\mu}_k$ in $S^{n-1}$ that approximates the directional measure $\mu$ of $Y$.

*Action*: Find a zonotope $\hat{Z}_k \in \mathcal{Z}^n$ that solves the following least squares problem:

$$\min_{Z \in \mathcal{Z}^n} \sum_{i=1}^{k} (y_i - h_Z(u_i))^2. \tag{77}$$

Calculate the finitely supported surface area measure $S(\hat{Q}_k, \cdot)$ of the origin-symmetric polytope $\hat{Q}_k$ satisfying

$$\hat{Z}_k = \Pi \hat{Q}_k$$

and set $\hat{\mu}_k = (1/2) S(\hat{Q}_k, \cdot)$. □

As was remarked for Algorithm NoisyBrightLSQ after the statement of that algorithm, $\hat{Z}_k$ can be assumed to be a sum of line segments, each parallel to a node corresponding to $U = \{u_1, \ldots, u_k\}$, and only the lengths of these line segments have to be determined. Applying the same observation leads to an output $\hat{\mu}_k$ that is supported by the finite set of nodes corresponding to $U$. This implementation of Algorithm NoisyRoseLSQ was suggested previously by Männle [23], who obtained the following result.

PROPOSITION 9.1. *Let $Y$ be a stationary fiber process in $\mathbb{R}^n$ with directional measure $\mu$ and let $(u_i)$ be a sequence in $S^{n-1}$ such that $(u_i^*)$ is evenly spread. If $\hat{\mu}_k$ is an output from Algorithm NoisyRoseLSQ as stated above, then, almost surely, $\hat{\mu}_k$ converges weakly to $\mu$, as $k \to \infty$.*

Männle [23] obtained Proposition 9.1 using local Kuhn–Tucker conditions for the solutions of a weighted least squares problem slightly more general than (77). However, the result follows immediately from Theorem 7.2 on observing that the map that takes $K \in \mathcal{K}^n$ to $S(K, \cdot)$ is weakly continuous on $\mathcal{K}^n$ (see, e.g., [28], page 205).

The remainder of this section is devoted to presenting a refinement of Proposition 9.1 that provides rates of convergence of the estimators. This requires the introduction of metrics on the cone of finite Borel measures in $S^{n-1}$ to quantify the deviation of the estimator from the true directional measure. Details for the following definitions in the case of probability measures can be found in Section 11.3 of [7]; the extension to arbitrary (nonnegative) measures is not difficult.

Let $\mu$ and $\nu$ be finite Borel measures in $S^{n-1}$. Define

$$(78) \qquad d_{\mathrm{D}}(\mu, \nu) = \sup\left\{ \left| \int_{S^{n-1}} f \, d(\mu - \nu) \right| : \|f\|_{BL} \leq 1 \right\},$$

where, for any real-valued function $f$ on $S^{n-1}$, we define

$$\|f\|_{BL} = \|f\|_\infty + \|f\|_L \quad \text{and} \quad \|f\|_L = \sup_{u \neq v} \frac{|f(u) - f(v)|}{\|u - v\|}.$$

It can be shown that $d_{\mathrm{D}}$ is a metric, sometimes called the *Dudley metric* (though he attributes its definition to Fortet and Mourier [9]) on the cone of finite Borel measures, inducing the weak topology. Now define

$$(79) \qquad \begin{aligned} d_{\mathrm{P}}(\mu, \nu) = \inf\{\varepsilon > 0 : &\mu(F) \leq \nu(F^\varepsilon) + \varepsilon, \\ &\nu(F) \leq \mu(F^\varepsilon) + \varepsilon, F \text{ closed in } S^{n-1}\}, \end{aligned}$$

where

$$F^\varepsilon = \left\{ u \in S^{n-1} : \inf_{v \in F} \|u - v\| < \varepsilon \right\}.$$

Then $d_\mathrm{P}$ is also a metric, the *Prohorov metric*, that induces the weak topology. The Dudley and Prohorov metrics are related, as we show below in Lemma 9.5.

The following proposition follows from a stability result of Hug and Schneider [17] that generalizes one step in the proof of the version of Proposition 2.2 due to Bourgain and Lindenstrauss [1].

PROPOSITION 9.2. *Let $K$ and $L$ be origin-symmetric convex bodies in $\mathbb{R}^n$, such that*

$$r_0 B \subset K, L \subset R_0 B,$$

*for some $0 < r_0 \leq R_0$. If $0 < b < 2/(n(n+4))$, there is a constant $c' = c'(b, n, r_0, R_0)$ such that*

(80) $$d_\mathrm{D}(S(K, \cdot), S(L, \cdot)) \leq c' \delta_2(\Pi K, \Pi L)^b.$$

PROOF. We refer the reader to Theorem 5.1 of [17]. In that result, more general than the statement of our theorem, take $\mu = S(K, \cdot) - S(L, \cdot)$ and $\Phi(u \cdot v) = |u \cdot v|$, so that according to [17], equation (52),

$$(T_\Phi(\mu))(u) = V(K|u^\perp) - V(L|u^\perp) = h_{\Pi K}(u) - h_{\Pi L}(u).$$

As is noted by Hug and Schneider [17], who assume throughout that $n \geq 3$, we may then take $\beta = (n+2)/2$ in their Theorem 5.1. With these substitutions, our theorem for $n \geq 3$ follows immediately.

When $n = 2$, Theorem 5.1 of [17] is still valid (and our theorem follows as before), but its proof requires an adjustment. One of the main steps is the approximation of a continuous function $f$ by its Poisson integral

$$f_r(u) = \frac{1}{V_{n-1}(S^{n-1})} \int_{S^{n-1}} \frac{1 - r^2}{(1 + r^2 - 2r u \cdot v)^{n/2}} f(v) \, dv,$$

where $0 < r < 1$ is a parameter. The proof of Theorem 5.1 of [17] uses the estimate

(81) $$\|f - f_r\|_\infty \leq 2^{n+1} \frac{V_{n-2}(S^{n-2})}{V_{n-1}(S^{n-1})} \|f\|_L (1 - r) \log \frac{2}{1 - r},$$

for $1/4 \leq r < 1$, from Lemma 5.5.8 of [15], where the proof applies only when $n \geq 3$. However, when $n = 2$ it can be shown that

(82) $$\|f - f_r\|_\infty \leq \frac{16\sqrt{3}}{\pi} \|f\|_L (1 - r) \log \frac{2}{1 - r}$$

for $1/4 \leq r < 1$. Although this estimate is slightly weaker than (81), it is sufficient to prove Theorem 5.1 of [17] for $n = 2$. For a proof of (82), see the Appendix of [11]. □

Let $\mathcal{D}$ denote the set of *degenerate* finite Borel measures in $S^{n-1}$, that is, those whose support is contained in a great sphere.

LEMMA 9.3. *Let $\mu$ be a finite Borel measure in $S^{n-1}$ and let*

(83) $$d_{\mathrm{D}}(\mu, \mathcal{D}) = \inf_{\nu \in \mathcal{D}} d_{\mathrm{D}}(\mu, \nu).$$

*Then the infimum is attained and the mapping $\mu \mapsto d_{\mathrm{D}}(\mu, \mathcal{D})$ is weakly continuous. Consequently, the support of $\mu$ is not contained in any great sphere if and only if $d_{\mathrm{D}}(\mu, \mathcal{D}) > 0$.*

PROOF. For $a \geq 0$, let
$$\mathcal{D}_a = \{\nu \in \mathcal{D} : \nu(S^{n-1}) \leq 2a\}.$$

If 0 denotes the zero measure, we have
$$d_{\mathrm{D}}(\mu, \mathcal{D}) \leq d_{\mathrm{D}}(\mu, 0) = \mu(S^{n-1}).$$

Therefore, if $a \geq \mu(S^{n-1})$, then

(84) $$d_{\mathrm{D}}(\mu, \mathcal{D}) = \inf_{\nu \in \mathcal{D}, d_{\mathrm{D}}(\mu,\nu) \leq a} d_{\mathrm{D}}(\mu, \nu) = \inf_{\nu \in \mathcal{D}_a} d_{\mathrm{D}}(\mu, \nu),$$

where the last equality comes from substituting $f \equiv 1$ in the definition (78) of $d_{\mathrm{D}}(\mu, \nu)$. It is easy to see that $\mathcal{D}$ is weakly closed and, hence, $\mathcal{D}_a$ is weakly compact, so the last infimum in (84) is attained.

Let $(\mu_k)$ be a sequence of finite Borel measures in $S^{n-1}$ converging to $\mu$. Choose $a$ so that $\mu_k(S^{n-1}) \leq a$ for all $k$. We know that there are measures $\nu \in \mathcal{D}_a$ and $\nu_k \in \mathcal{D}_a$, $k = 1, 2, \ldots$, such that
$$d_{\mathrm{D}}(\mu, \mathcal{D}) = d_{\mathrm{D}}(\mu, \nu) \quad \text{and} \quad d_{\mathrm{D}}(\mu_k, \mathcal{D}) = d_{\mathrm{D}}(\mu_k, \nu_k),$$

for $k = 1, 2 \ldots$. The weak compactness of $\mathcal{D}_a$ implies that a subsequence of $(\nu_k)$ converges to a measure $\tilde{\nu} \in \mathcal{D}_a$. Then
$$\begin{aligned}
d_{\mathrm{D}}(\mu, \mathcal{D}) &\leq d_{\mathrm{D}}(\mu, \tilde{\nu}) \\
&= \liminf_{k \to \infty} d_{\mathrm{D}}(\mu_k, \nu_k) \\
&= \liminf_{k \to \infty} d_{\mathrm{D}}(\mu_k, \mathcal{D}) \\
&\leq \limsup_{k \to \infty} d_{\mathrm{D}}(\mu_k, \mathcal{D}) \\
&\leq \limsup_{k \to \infty} d_{\mathrm{D}}(\mu_k, \nu) \\
&= d_{\mathrm{D}}(\mu, \nu) = d_{\mathrm{D}}(\mu, \mathcal{D}).
\end{aligned}$$

Therefore,
$$\lim_{k \to \infty} d_{\mathrm{D}}(\mu_k, \mathcal{D}) = d_{\mathrm{D}}(\mu, \mathcal{D}),$$

as required. □

The following refinement of Proposition 9.1 is phrased in terms of the Dudley metric. For $n \geq 3$, the extra condition that $d \leq d_{\mathrm{D}}(\mu, \mathcal{D})$ for some $d > 0$ is needed. It is a natural analog of the condition that $rB \subset K$ for some $r > 0$ in earlier results, such as Theorem 7.6. Lemma 9.3 implies that such a lower bound $d > 0$ always exists due to our general assumption that the directional measure $\mu$ is not degenerate.

THEOREM 9.4. *Let $\sigma > 0$. Let $Y$ be a stationary fiber process in $\mathbb{R}^n$ with directional measure $\mu$ and length density $\overline{L} = \mu(S^{n-1})$. Let $(u_i)$ be a sequence of directions in $S^{n-1}$ with $\Delta_k^* = O(k^{-1/(n-1)})$ and let $\hat{\mu}_k$ be an output from [Algorithm NoisyRoseLSQ](#) as stated above.*

*If $n = 2$ and $\beta > 0$, then, almost surely, there are constants $C_{14} = C_{14}(\sigma, \overline{L}, (u_i), \beta)$ and $N_{14} = N_{14}(\sigma, \overline{L}, (u_i), \beta)$ such that*

$$d_{\mathrm{D}}(\mu, \hat{\mu}_k) \leq C_{14} k^{-2/15 + \beta}, \tag{85}$$

*for $k \geq N_{14}$.*

*For $n \geq 3$, let $0 < d \leq d_{\mathrm{D}}(\mu, \mathcal{D})$. If $n = 3$ or $4$ and $\beta > 0$, then, almost surely, there are constants $C_{15} = C_{15}(\sigma, n, \overline{L}, d, (u_i), \beta)$ and $N_{15} = N_{15}(\sigma, n, \overline{L}, d, (u_i), \beta)$ such that*

$$d_{\mathrm{D}}(\mu, \hat{\mu}_k) \leq C_{15} k^{-(n+2)/((n+4)(2n+1)) + \beta}, \tag{86}$$

*for $k \geq N_{15}$.*

*Finally, if $n \geq 5$ and $\beta > 0$, there are constants $C_{16} = C_{16}(\sigma, n, \overline{L}, d, (u_i), \beta)$ and $N_{16} = N_{16}(\sigma, n, \overline{L}, d, (u_i), \beta)$ such that*

$$d_{\mathrm{D}}(\mu, \hat{\mu}_k) \leq C_{16} k^{-2/((n-1)(n+4)) + \beta}, \tag{87}$$

*for $k \geq N_{16}$.*

PROOF. Let $K$ and $\hat{Q}_k$ be the origin-symmetric convex bodies with surface area measures $2\mu$ and $2\hat{\mu}_k$, respectively, and recall that $\hat{Z}_k = \Pi \hat{Q}_k$.

Suppose that $n = 2$. According to [28], pages 290–291, the mean width

$$\overline{w}(K) = \frac{1}{\pi} \int_{S^1} h_K(u) \, du$$

of $K$ satisfies $\pi \overline{w}(K) = S(K, S^1) = \overline{L}$. Since $K = -K$, for each $x \in K$, $[-x, x] \subset K$ and so

$$\frac{4}{\pi} \|x\| = \overline{w}([-x, x]) \leq \overline{w}(K) = \frac{1}{\pi} \overline{L}.$$

It follows that $K \subset (\overline{L}/4)B$. By the case $n = 2$ of Theorem 7.5, with

$$R = \max\{\overline{L}/4, (2^{15/2} \sigma / \kappa_{n-1})^{1/(n-1)}\},$$

almost surely, there are constants $C_{17} = C_{17}(\sigma, \overline{L}, (u_i))$ and $N_{17} = N_{17}(\sigma, \overline{L}, (u_i))$ such that

$$\delta_2(\Pi K, \Pi \hat{Q}_k) = \delta_2(\Pi K, \hat{Z}_k) \leq C_{17} k^{-2/5},$$

for all $k \geq N_{17}$. Inequality (80) with $L = \hat{Q}_k$ now implies (85).

Suppose that $n \geq 3$ and let $0 < d \leq d_D(\mu, \mathcal{D})$, which is possible by Lemma 9.3. Let $\mathcal{M} = \mathcal{M}(\overline{L}, d)$ be the set of all finite Borel measures $\nu$ in $S^{n-1}$ such that $\nu(S^{n-1}) \leq \overline{L}$ and $d_D(\nu, \mathcal{D}) \geq d$. Then $\mu \in \mathcal{M}$ and $\mathcal{M}$ is weakly compact by Lemma 9.3. Using the equicontinuity of the family $\{f_u : u \in S^{n-1}\}$ of functions defined by $f_u(v) = |u \cdot v|$ for $v \in S^{n-1}$, we see that the map $T: S^{n-1} \times \mathcal{M} \to \mathbb{R}$ defined by

$$T(u, \nu) = \int_{S^{n-1}} |u \cdot v| \, d\nu(v)$$

is continuous. Therefore, $T$ attains its minimum $r = r(n, \overline{L}, d)$ at some point $(u_0, \nu_0)$ in the compact set $S^{n-1} \times \mathcal{M}$. Note that

$$r = T(u_0, \nu_0) = \int_{S^{n-1}} |u_0 \cdot v| \, d\nu_0(v) > 0,$$

as $\nu_0$ is not degenerate. Then $T(u, \mu) \geq r$ for all $u \in S^{n-1}$, so by (75) and the fact that $S(K, \cdot) = 2\mu$, we have

$$h_{\Pi K}(u) = \int_{S^{n-1}} |u \cdot v| \, d\mu(v) \geq r,$$

for all $u \in S^{n-1}$. Therefore, $rB \subset \Pi K$. On the other hand, (75) also implies $h_{\Pi K}(u) \leq \overline{L} = R$. Summarizing, we have shown that there are constants $0 < r < R$, depending only on $n$, $\overline{L}$ and $d$, such that

$$rB \subset \Pi K \subset RB.$$

As in the proof of Theorem 7.2, we can conclude

$$r_0 B \subset K \subset R_0 B,$$

with positive constants $r_0 = r_0(n, d, R)$ and $R_0 = R_0(n, d, R)$. Theorem 7.5 can now be applied with

$$R = \max\{R_0, (2^{15/2} \sigma / \kappa_{n-1})^{1/(n-1)}\}.$$

Inequality (80) with $L = \hat{Q}_k$ then yields (86) and (87). $\square$

To obtain a version of Theorem 9.4 in terms of the Prohorov metric, the following lemma is useful.

LEMMA 9.5. *Let $\mu$ and $\nu$ be finite Borel measures in $S^{n-1}$ with $m_0 = \mu(S^{n-1}) \neq 0$. If $d_D(\mu, \nu) \leq 1$, then*

$$d_P(\mu, \nu) \leq (1 + \sqrt{3 + m_0}) \, d_D(\mu, \nu)^{1/2}.$$

PROOF. We may assume that $n_0 = \nu(S^{n-1}) \neq 0$ and let $\mu_1 = \mu/m_0$ and $\nu_1 = \nu/n_0$. For $s, t > 0$, the definition (79) of $d_P$ implies that

$$d_P(s\nu, t\nu) = \inf\{\varepsilon > 0 : s\nu(F) - t\nu(F^\varepsilon) \leq \varepsilon, t\nu(F) - s\nu(F^\varepsilon) \leq \varepsilon,$$
$$F \text{ closed in } S^{n-1}\}$$
$$\leq \inf\{\varepsilon > 0 : s\nu(F^\varepsilon) - t\nu(F^\varepsilon) \leq \varepsilon, t\nu(F^\varepsilon) - s\nu(F^\varepsilon) \leq \varepsilon,$$
$$F \text{ closed in } S^{n-1}\}$$
$$\leq n_0 |s - t|,$$

while the definition (78) of $d_D$ (with $f \equiv 1$) yields

$$d_D(\mu, \nu) \geq |m_0 - n_0|.$$

Therefore,

(88)
$$\begin{aligned} d_P(\mu, \nu) &\leq d_P(\mu, (m_0/n_0)\nu) + d_P((m_0/n_0)\nu, \nu) \\ &\leq d_P(m_0\mu_1, m_0\nu_1) + n_0 |m_0/n_0 - 1| \\ &\leq d_P(m_0\mu_1, m_0\nu_1) + d_D(\mu, \nu). \end{aligned}$$

Let $0 < \varepsilon < d_P(m_0\mu_1, m_0\nu_1)$. By (79), there is a closed set $F$ in $S^{n-1}$ such that

$$\mu_1(F) > \nu_1(F^\varepsilon) + \frac{\varepsilon}{m_0} \quad \text{or} \quad \nu_1(F) > \mu_1(F^\varepsilon) + \frac{\varepsilon}{m_0}.$$

Setting $\alpha = \varepsilon/m_0$ and $\beta = \varepsilon$ in Proposition 3 of [6], we obtain

(89)
$$\frac{2}{m_0(2+\varepsilon)}\varepsilon^2 \leq d_D(\mu_1, \nu_1).$$

By (79) again, we have $d_P(\mu_1, \nu_1) \leq 1$ and

$$d_P(m_0\mu_1, m_0\nu_1) = \inf\left\{\varepsilon > 0 : \mu_1(F) \leq \nu_1(F^\varepsilon) + \frac{\varepsilon}{m_0},\right.$$
$$\left.\nu_1(F) \leq \mu_1(F^\varepsilon) + \frac{\varepsilon}{m_0}, F \text{ closed in } S^{n-1}\right\}.$$

Therefore, if $m_0 \leq 1$, we have

$$d_P(m_0\mu_1, m_0\nu_1) \leq d_P(\mu_1, \nu_1) \leq 1,$$

while, if $m_0 \geq 1$, then

$$d_P(m_0\mu_1, m_0\nu_1)$$
$$= m_0 \inf\{\varepsilon > 0 : \mu_1(F) \leq \nu_1(F^{m_0\varepsilon}) + \varepsilon,$$
$$\nu_1(F) \leq \mu_1(F^{m_0\varepsilon}) + \varepsilon, F \text{ closed in } S^{n-1}\}$$
$$\leq m_0 d_P(\mu_1, \nu_1) \leq m_0.$$

Thus, for any $m_0 > 0$, we have $\varepsilon < d_{\mathrm{P}}(m_0\mu_1, m_0\nu_1) \leq 1 + m_0$. Substitution into (89) yields

$$\frac{2}{m_0(3+m_0)}\varepsilon^2 \leq d_{\mathrm{D}}(\mu_1, \nu_1).$$

As $\varepsilon < d_{\mathrm{P}}(m_0\mu_1, m_0\nu_1)$ was arbitrary, we conclude that

$$\frac{2}{m_0(3+m_0)}d_{\mathrm{P}}(m_0\mu_1, m_0\nu_1)^2 \leq d_{\mathrm{D}}(\mu_1, \nu_1)$$

$$\leq d_{\mathrm{D}}(\mu/m_0, \nu/m_0) + d_{\mathrm{D}}(\nu/m_0, \nu/n_0)$$

$$\leq \frac{1}{m_0}d_{\mathrm{D}}(\mu, \nu) + n_0\left|\frac{1}{m_0} - \frac{1}{n_0}\right|$$

$$\leq \frac{2}{m_0}d_{\mathrm{D}}(\mu, \nu).$$

Substituting this into (88), and using the hypothesis $d_{\mathrm{D}}(\mu, \nu) \leq 1$, we obtain the desired inequality. $\square$

With Lemma 9.5 in hand, the estimates of Theorem 9.4 can be converted to the Prohorov metric. Since this is routine, we shall only give one example. By Lemma 9.5 and (86), under the hypotheses of Theorem 9.4 with $n=3$, for all $\gamma > 0$, almost surely, there are constants $C_{18} = C_{18}(\sigma, R, (u_i), \gamma)$ and $N_{18} = N_{18}(\sigma, R, (u_i), \gamma)$ such that

$$d_{\mathrm{P}}(\mu, \hat{\mu}_k) \leq C_{18} k^{-5/98+\gamma},$$

for all $k \geq N_{18}$. Thus, the exponent when $n=3$ is approximately $-1/20$.

We close this section with a comment on the assumption in (76) that the errors are normally distributed. In applications, the measurements $y_i$ come from counting intersection points, so they are integer random variables. If $\overline{L}$ is large, our assumption is appropriate. Otherwise, a model that allows only integer values for $y_i$ could be more apt. For example, if $Y$ is a Poisson line process (one of the most common models in stochastic geometry), then the number $y_i$ of intersection points of its fibers with a unit window in $u^\perp$ is Poisson distributed with mean $\gamma(u_i)$, $i = 1, \ldots, k$. Under this assumption on the distribution, the maximum likelihood problem no longer corresponds to a quadratic program. Nevertheless, its solution is a strongly consistent estimator for $\mu$ (see [20]), and the tools provided by van de Geer [32] would still allow results giving rates of convergence.

**Acknowledgments.** We thank Chris Eastman, Amyn Poonawala, Greg Richardson, Thomas Riehle and Chris Street for writing the computer programs for Monte Carlo simulations. We are also grateful to Sara van de Geer for a helpful correspondence concerning her work and to the referees for comments that led to a better presentation of our results.

R. J. GARDNER
DEPARTMENT OF MATHEMATICS
WESTERN WASHINGTON UNIVERSITY
BELLINGHAM, WASHINGTON 98225-9063
USA
E-MAIL: Richard.Gardner@wwu.edu

M. KIDERLEN
DEPARTMENT OF MATHEMATICAL SCIENCES
UNIVERSITY OF AARHUS
NY MUNKEGADE, DK-8000 AARHUS C
DENMARK
E-MAIL: kiderlen@imf.au.dk

P. MILANFAR
DEPARTMENT OF ELECTRICAL ENGINEERING
UNIVERSITY OF CALIFORNIA
1156 HIGH STREET
SANTA CRUZ, CALIFORNIA 95064-1077
USA
E-MAIL: milanfar@soe.ucsc.edu